\documentclass{amsart}
\usepackage{amssymb}
\usepackage{graphics}
\usepackage{color}

\def\2{\color{red}}
\let\3\relax

\usepackage[hyperref]{degt}

\makeatletter
\def\HyPsd@CatcodeWarning#1{}
\makeatother

\def\CG#1{\Z/#1}

\let\Gm\mu
\let\Gn\nu
\let\Gu\upsilon
\let\oL\Lambda
\let\oO\Omega
\def\oM{\mathrm{M}}
\def\oN{\mathrm{N}}
\def\oH{\mathrm{H}}
\def\oE{\bar{\mathrm{E}}}
\def\ooE{\mathrm{E}}
\def\ooO{\mathrm{O}}

\def\oZ{\oM_0}

\def\e{\mathrm{e}}
\def\be{\bar\e}
\def\rt{\frak{t}}

\def\bX{\bold X}
\def\bQ{\bold Q}

\def\bH{\bold H}
\def\bS{\bold S}
\def\bL{\bold L}
\def\GQ{\operatorname{GQ}}
\def\kchar{\fchar\Bbbk}
\def\CF{\Cal F}

\def\F{\Bbb F}
\def\pencil{\Cal P}
\def\val{\operatorname{val}}
\def\mult{\operatorname{mult}}

\def\Fn{\operatorname{Fn}}
\def\fb{\operatorname{fb}}
\def\fib{\operatorname{fib}}
\def\dif{\mathrm{d}}

\def\aD{\smash{\tilde D}}
\def\tP{\smash{\tilde P}}
\def\tS{\smash{\tilde S}}
\def\CK{\Cal K}
\def\CL{\Cal L}
\def\CS{\Cal S}
\def\CT{\Cal T}
\def\CX{\Cal X}
\def\CA{\Cal A}
\def\CH{\Cal H}
\def\DD{\Cal D}
\def\CE{\Cal E}
\def\CU{\Cal U}
\def\CV{\Cal V}
\def\CO{\Cal O}
\def\barT{\bar T}

\def\bbE{\bar\bE}

\def\qc{q_\circ}
\def\cc{\circ}

\def\Br{\operatorname{Br}}
\def\dual{^\vee}

\let\pv=\epsilon

\makeatletter
\def\PS#1{\expandafter\@PSi#1\end@PS}
\def\@PSi#1[#2[#3]#4,#5]{(#3)^{#5}\futurelet\next\@PSii}
\def\@PSii{\ifx\next,\expandafter\@PSi\fi}
\let\end@PS\relax
\makeatother
\def\ps{\frak{p}}
\let\Artin\sigma
\def\Fano{\Cal F}

\def\sect#1{\QOPNAME{s}(#1)}
\def\flines#1{\QOPNAME{ln}(#1)}

\begin{document}


\title{Lines in supersingular quartics}

\author{Alex Degtyarev}

\address{%
Department of Mathematics\\
Bilkent University\\
06800 Ankara, TURKEY}

\email{
degt@fen.bilkent.edu.tr}

\thanks{%
The author was supported by the JSPS grant L15517
and T\"{U}B\DOTaccent{I}TAK grant 114F325}

%
%



\begin{abstract}
We show that the number of lines contained in a supersingular quartic surface
is 40 or at most 32, if the characteristic of the field equals 2, and it is
112, 58, or at most 52, if the characteristic equals 3.
If the quartic is not supersingular, the number of lines is at most 60
in both cases.
We also give a complete classification of large configurations of lines.
\end{abstract}

\maketitle

\section{Introduction}\label{S.intro}

Throughout the paper,
unless specified otherwise,
$X$ stands for a nonsingular quartic surface in
the projective space $\Cp3$
over an algebraically closed field~$\Bbbk$.

\subsection{Motivation}\label{s.motivation}
A simple dimension count shows that, unlike quadrics or cubics, a generic
quartic surface $X\subset\Cp3$ contains no straight lines.
On the other hand, it has been known since F.~Schur~\cite{Schur:quartics} that
there exists a quartic~$X_{64}$ containing $64$ lines.
B.~Segre~\cite{Segre} proved that the number~$64$ is maximal possible.
After a period of oblivion, S.~Rams and M.~Sch\"{u}tt~\cite{rams.schuett}
bridged a gap in Segre's arguments and extended his (correct) bound~$64$ to
any algebraically closed field of characteristic $\kchar\ne2,3$.
Since Schur's quartic~$X_{64}$ has a nonsingular reduction over such fields,
the bound is sharp. If $\kchar=3$, the maximal number of lines is~$112$,
see~\cite{rams.schuett:char3}; if $\kchar=2$, the maximal number is~$60$, see
\autoref{th.ordinary} below.

At the same time, in recent paper~\cite{DIS}, we suggested an alternative
approach to Segre's theorem over~$\C$, using the theory of $K3$-surfaces
and Nikulin's theory of discriminant forms~\cite{Nikulin:forms}.
We reestablished Segre's bound~$64$, proved that Schur's quartic~$X_{64}$ is
the only one containing $64$ lines (see \autoref{cor.Schur}
below for a similar
statement over an arbitrary field), and gave a complete classification of all
large configurations: up to projective equivalence, there are but
ten quartics containing more than $52$ lines. (Other results of~\cite{DIS}
are the sharp bound~$56$ for the number of \emph{real} lines in a
\emph{real} quartic and the bound $52$ for the number of lines defined
over~$\Q$ in a quartic defined over~$\Q$.)

In the present paper, we obtain similar refined results for the cases
$\kchar=2$ or~$3$. According to~\cite{Lieblich.Maulik}, if a quartic~$X$ is
\emph{not} supersingular, it is subject to
the same
lattice theoretical restrictions as quartics defined over~$\C$.
Hence, the list of large configurations found in~\cite{DIS} applies to such
quartics as a ``bound'', with some entries missing over some fields. (An example of such
missing entries is \autoref{th.ordinary}, which rules out the Schur
configuration~$\bX_{64}$ in characteristics~$2$ and~$3$.)
Therefore, we concentrate on supersingular surfaces; our principal
results, Theorems~\ref{th.char=2} and~\ref{th.char=3}, show that the
configurations of lines realized by
such surfaces
do differ dramatically from
the eight configurations found in~\cite{DIS}.

Characteristics~$2$ and~$3$ are naturally special for quartics:
these primes divide
the degrees of the defining polynomial and its derivatives, and it is these
(and only these)
characteristics where pencils of curves of arithmetic genus~$1$ ---one of the
principal tools commonly used in the theory--- may
become quasi-elliptic.
Note though that there also are interesting supersingular quartics over
other fields: thus, the quartic in characteristic~$7$ discussed in
\autoref{rem.triangle.free} beats by~$10$ all other known examples of
the so-called triangle free configurations.
Phenomena specific to fields of other characteristics
will be the subject of a forthcoming paper.

\subsection{Principal results}\label{s.results}

The set of lines in a quartic~$X$ is denoted by $\Fn X$, and the sublattice
spanned by the classes of the lines and plane section is denoted by
$\Fano(X)\subset\NS(X)$. We denote by $\Artin:=\Artin(X)$ the Artin invariant
of a supersingular $K3$-surface~$X$ (see \autoref{th.NS}).
An important easily comparable
combinatorial invariant of a configuration of lines is its
\emph{pencil structure}~$\ps$, \ie, the list of the types $(p,q)$ of
all pencils
$\pencil(l)$,
$l\in\Fn X$ (see \autoref{s.pencils}).
We use the partition notation, a
``factor'' $(p,q)^m$ standing for $m$ copies of the type $(p,q)$.

In the statements, we identify ``interesting'' quartics by the triple
$(\ps,\Artin,\rank\Fano)$:
these triples suffice to distinguish all examples found in the paper.
More details, such as the Gram matrix of
$\Fano(X)$ and coordinates of the lines in the
N\'{e}ron--Severi
lattice
$\NS(X)$, are available from the author in electronic form.

The principal results of the paper are \autoref{th.char=2} (supersingular
quartics in characteristic~$2$) and \autoref{th.char=3} (supersingular
quartics in characteristic~$3$). In \autoref{th.ordinary},
we reduce Segre's bound for quartics that are not supersingular.

\theorem[see \autoref{proof.char=2}]\label{th.char=2}
Assume that $\kchar=2$ and $X$ is supersingular.
Then either $\ls|\Fn X|=40$, and
there are at most five configurations\rom:
\roster
\item\label{40.1}
$\ps=\PS{ [ [ 2, 6 ], 40 ] }$, $\Artin(X)=3$, $\rank\Fano(X)=22$,
\item\label{40.2}
$\ps=\PS{ [ [ 2, 6 ], 40 ] }$, $\Artin(X)=3$, $\rank\Fano(X)=21$,
\item\label{40.3}
$\ps=\PS{ [ [ 4, 0 ], 4 ], [ [ 2, 6 ], 36 ] }$, $\Artin(X)=3$, $\rank\Fano(X)=22$,
\item\label{40.4}
$\ps=\PS{ [ [ 4, 0 ], 8 ], [ [ 2, 6 ], 32 ] }$, $\Artin(X)=3$, $\rank\Fano(X)=20$,
\item\label{40.5}
$\ps=\PS{ [ [ 4, 0 ], 40 ] }$, $\Artin(X)=3$, $\rank\Fano(X)=16$,
\endroster
or $\ls|\Fn X|\le32$.
\endtheorem

Due
to the lack of the existence statement in \autoref{th.Saint-Donat}
(\cf. \autoref{rem.Saint-Donat} below), we do
not assert that the five configurations listed in \autoref{th.char=2} are
realizable: one would have to find explicit defining equations for these
surfaces, but this task is beyond the scope of the present paper.
Conjecturally, the last bound $\ls|\Fn X|\le32$ is sharp: there are examples
of configurations with $32$ lines, but their realizability
by smooth supersingular quartics is also open.

\theorem[see \autoref{proof.char=3}]\label{th.char=3}
Assume that  $\kchar=3$ and $X$ is supersingular.
Then either $\ls|\Fn X|=112$, and $X$ is the Fermat quartic\rom:
\roster
\item\label{112.(10,0)}
$\ps=\PS{ [ [ 10, 0 ], 112 ]}$, $\Artin(X)=1$, $\rank\Fano(X)=22$,
\endroster
or $\ls|\Fn X|=58$, and there are three configurations\rom:
\roster[2]
\item\label{58.(10,0)x2}
$\ps=\PS{[ [ 10, 0 ], 2 ], [ [ 1, 9 ], 54 ], [ [ 1, 0 ], 2 ]}$,
$\Artin(X)=2$, $\rank\Fano(X)=22$,
\item\label{58.(10,0)}
$\ps=\PS{[ [ 10, 0 ], 1 ], [ [ 4, 6 ], 27 ], [ [ 4, 0 ], 12 ], [ [ 1, 9 ], 18 ]}$,
$\Artin(X)=2$, $\rank\Fano(X)=22$,
\item\label{58.(7,0)}
$\ps=\PS{ [ [ 7, 0 ], 2 ], [ [ 4, 6 ], 18 ], [ [ 3, 6 ], 36 ], [ [ 1, 9 ], 2 ]}$,
$\Artin(X)=2$, $\rank\Fano(X)=21$,
\endroster
or $\ls|\Fn X|\le52$, and this bound is sharp.
\endtheorem

D.~Veniani (private communication) has found explicit defining equations of
the three quartics with $58$ lines. Alternatively, quartics as in
\autoref{th.char=3}\iref{58.(10,0)x2} and~\iref{58.(10,0)} are described in
Propositions~\ref{prop.(10,0)++} and~\ref{prop.(10,0)+(4,6)},
respectively.

In \autoref{th.char=2}\iref{40.5},
the configuration $\Fn X$ constitutes (in the sense
described in \autoref{s.GQ} below)
the so-called
generalized quadrangle $W(3)$.
In \autoref{th.char=3}\iref{112.(10,0)}, $\Fn X$ constitutes the only
generalized quadrangle $\GQ(3,9)\cong Q(5,3)$.

If $X$ is not supersingular, the situation also differs from that in
characteristic~$0$, as some quartics defined over algebraic number fields
become singular and/or acquire extra lines when reduced to positive
characteristics. We have the following bound;
its sharpness is discussed in \autoref{rem.ordinary}.

\theorem[see \autoref{proof.ordinary}]\label{th.ordinary}
Assume that  $\kchar=2$ or~$3$ and $X$ is not supersingular.
Then $\ls|\Fn X|\le60$.
\endtheorem


According to~\cite{DIS}, there are considerable gaps
in the set of values taken by the number of lines in a nonsingular quartic
defined over~$\C$.
We conjecture similar gaps
for supersingular quartics
in characteristics~$2$ and~$3$.

\conjecture[see \autoref{rem.char=2.counts}]\label{conj.char=2}
Under the assumptions of \autoref{th.char=2},
the number $\ls|\Fn X|$ takes values in
the set
$\{0,1,\ldots,17,18, 20, 22, 24, 28, 32, 40\}$.
\endconjecture

\conjecture[see \autoref{rem.char=3.counts}]\label{conj.char=3}
Under the assumptions of \autoref{th.char=3},
one has $\ls|\Fn X|=1\bmod3$ whenever $\ls|\Fn X|\ge40$.
\endconjecture

\subsection{Contents of the paper}
Sections~\ref{S.lattices} and~\ref{S.K3} are preliminary: we summarize the
necessary facts concerning integral lattices, discriminant forms,
$K3$-surfaces, and (quasi-)elliptic pencils.
In \autoref{S.config}, we summarize and extend some intermediate results
of~\cite{DIS}, introducing the principal technical tools---configurations and
pencils. Then,
in \autoref{S.trig.free}, we treat the so-called triangle free
configurations, also following~\cite{DIS}. We prove a (rather week)
characteristic independent bound and a few intermediate lemmas that are used
later.
The principal results of the paper, \viz. Theorems~\ref{th.char=3}
and~\ref{th.char=2}, are proved in \autoref{S.exotic} and
\autoref{S.2-exotic}, where we study in detail pencils in supersingular
quartics over fields of characteristic~$3$ and~$2$, respectively.
Finally, in \autoref{S.geometry}, intuitive geometric arguments are
used to rule out Schur's configuration~$\bX_{64}$
in characteristics~$3$ and~$2$ and prove
\autoref{th.ordinary}; we conclude this section with explicit defining
equations of several supersingular quartics in characteristic~$3$.

\subsection{Acknowledgements}
I cordially thank
Dmitrii Pasechnik, Matthias Sch\"{u}tt, Tetsuji Shioda, and Davide Veniani
for a number of comments, suggestions, and fruitful and motivating
discussions.
My special gratitude goes to Ichiro Shimada, who introduced me to the world
of supersingular $K3$-surfaces and generously shared his ideas concerning
this project.
I thank
the anonymous referee of this paper, who pointed out to a few
inaccuracies in the text.
This paper was written during my sabbatical stay at
Hiroshima University, supported by the Japan Society for the Promotion of
Science; I am grateful to these institutions for their hospitality and
support.

\section{Lattices}\label{S.lattices}

In this introductory section we recall briefly a few elementary facts
concerning integral lattices and their discriminant forms. The principal
reference is~\cite{Nikulin:forms}.

\subsection{Finite quadratic forms\pdfstr{}{ {\rm(see~\cite{Miranda.Morrison:book,Nikulin:forms})}}}\label{s.finite.forms}
A \emph{finite quadratic form} is a finite abelian group~$\CL$ equipped with
a map $q\:\CL\to\Q/2\Z$ quadratic in the sense that
\[*
q(x+y)=q(x)+q(y)+2b(x,y),\quad
q(nx)=n^2q(x),\quad
x,y\in\CL,\ n\in\Z,
\]
where $b\:\CL\otimes\CL\to\Q/\Z$ is a symmetric bilinear form and
$2\:\Q/\Z\to\Q/2\Z$ is the natural isomorphism.
We often abbreviate $x^2:=q(x)$ and $x\cdot y:=b(x,y)$. Clearly, $b$
is determined by~$q$; the converse holds if and only if $\ls|\CL|$ is prime
to~$2$.

There is a direct sum decomposition $\CL=\bigoplus_p\CL_p$, where
$\CL_p:=\CL\otimes\Z_p$ is the $p$-primary part of~$\CL$ and $p$ runs over
all primes. The \emph{length} $\ell(\CL)$ is the minimal number of generators
of~$\CL$; we abbreviate $\ell_p(\CL):=\ell(\CL_p)$.
The form on~$\CL_2$ is called \emph{even} if $x^2=0\bmod\Z$ for each element
$x\in\CL_2$ of order~$2$; otherwise, it is called \emph{odd}.
There is a unique vector $c\in\CL_2/2\CL_2$ with the property
$x^2=x\cdot c\bmod\Z$ for each element $x\in\CL_2$ or order~$2$; it is called
the \emph{characteristic vector}. The form on $\CL_2$ is even if and only if
$c=0$.

A finite quadratic/bilinear form is \emph{nondegenerate} if the
\emph{associated map}
\[*
\CL\longto\Hom(\CL,\Q/\Z),\quad
 x\longmapsto(y\mapsto x\cdot y)
\]
is an isomorphism.
A nondegenerate finite quadratic form splits into an orthogonal direct sum
of cyclic forms $\<\frac{m}{n}\>$, $\gcd(m,n)=1$, $mn=0\bmod2$
(defined on the cyclic group $\CG{m}$) and length~$2$ blocks
(on the group $(\CG{n})^2$)
\[*
\CU_n:=\biggl<\begin{matrix}0&1/n\\1/n&0\end{matrix}\biggr>,\quad
\CV_n:=\biggl<\begin{matrix}2/n&1/n\\1/n&2/n\end{matrix}\biggr>,\quad
\text{where $n=2^k$, $k\ge1$}.
\]

Given a prime~$p$,
the \emph{determinant} $\det_p\CL$ of a
nondegenerate
finite quadratic form
is the determinant of the matrix of the form on the $p$-group~$\CL_p$
in any minimal
basis. According to~\cite{Miranda.Morrison:book}, one has
$\det_p\CL=u\ls|\CL_p|\1$, where $u\in\Z_p^\times$; the unit~$u$ is well
defined modulo $(\Z_p^\times)^2$ unless $p=2$ and $\CL_2$ is odd; in the
latter case, $\det_2\CL$ is well defined modulo the subgroup generated by
$(\Z_2^\times)^2$ and~$5$.

The \emph{Brown invariant} of a nondegenerate finite quadratic form~$\CL$ is the
residue $\Br q=\Br\CL\in\CG8$ defined by the Gauss sum
\[*
\exp\bigl(\tfrac14i\pi\Br\CL\bigr)=
 \ls|\CL|^{-\frac12}\sum_{x\in\CL}\exp\bigl(i\pi x^2\bigr).
\]
The Brown invariant is additive: $\Br(\CL'\oplus\CL'')=\Br\CL'+\Br\CL''$.

A finite quadratic form~$q$ (respectively, bilinear form~$b$) on~$\CL$ is
\emph{null-cobordant} if there exists a $q$-isotropic (respectively,
$b$-isotropic) subgroup $\CK\subset\CL$ of maximal order, \ie,
such that $\ls|\CL|=\ls|\CK|^2$ or, equivalently, $\CK=\CK^\perp$.
If a quadratic form~$q$ is null-cobordant,
then $\Br q=0$; if $q$ is defined on a $2$- or $3$-elementary group, the
converse also holds.
More generally, for any $q$-isotropic subgroup $\CK\subset\CL$
one has the identity
$\Br(\CK^\perp\!/\CK)=\Br\CL$ (\cf. \autoref{th.extension} below).

\subsection{Integral lattices\pdfstr{}{ {\rm(see~\cite{Nikulin:forms})}}}\label{s.lattice}
An \emph{\rom(integral\rom) lattice} is a finitely generated free abelian
group~$L$ equipped with a symmetric bilinear form $b\:L\otimes L\to\Z$;
usually,
we abbreviate $x^2:=b(x,x)$ and $x\cdot y:=b(x,y)$.
A lattice~$L$ is \emph{even} if $x^2=0\bmod2$ for all $x\in L$; otherwise,
$L$ is \emph{odd}.
The \emph{determinant} $\det L\in\Z$ is the determinant of the matrix of~$b$
in any integral basis. The lattice~$L$ is called \emph{nondegenerate} if
$\det L\ne0$; it is called \emph{unimodular} if $\det L=\pm1$.
Equivalently, $L$ is nondegenerate if and only if its \emph{kernel}
\[*
\ker L=L^\perp:=\bigl\{x\in L\bigm|\text{$x\cdot y=0$ for all $y\in L$}\bigr\}
\]
is trivial.
A \emph{characteristic vector} of a unimodular lattice~$L$ is a vector
$u\in L$ such that $x^2=x\cdot u\bmod2$ for all $x\in L$. Such a vector
exists and is unique $\bmod\,2L$.

The inertia indices $\Gs_\pm$ and signature $\Gs:=\Gs_+-\Gs_-$
of a lattice~$L$ are
defined  as those of $L\otimes\Q$; a nondegenerate lattice~$L$ is called
\emph{hyperbolic} if $\Gs_+L=1$.

Let $L$ be a nondegenerate lattice. Then, we have a canonical inclusion
\[
L\subset L\dual:=\Hom(L,\Z)=
 \bigl\{x\in L\otimes\Q\bigm|\text{$x\cdot y\in\Z$ for all $y\in L$}\bigr\}.
\label{eq.dual}
\]
The finite group $\discr L:=L\dual\!/L$ of order $\ls|\det L|$ is called the
\emph{discriminant \rom(group\rom)} of~$L$.
This group inherits from
$L\otimes\Q$ the symmetric bilinear \emph{discriminant form}
\[*
b\:(x\bmod L)\otimes(y\bmod L)\mapsto(x\cdot y)\bmod\Z\in\Q/\Z
\]
and, if $L$ is even, its quadratic extension
\[*
q\:(x\bmod L)\mapsto x^2\bmod2\Z\in\Q/2\Z.
\]
These forms are taken into account whenever we speak about (anti-)isometries
of discriminant groups.
We abbreviate $\discr_pL:=(\discr L)\otimes\Z_p$.
A lattice~$L$ is said to be \emph{$p$-elementary} if $\discr L$ is a
$p$-elementary group.

To avoid confusion, we fix the notation:
\roster*
\item
$L^n$, $n\in\N$, is the orthogonal direct sum of $n$ copies of~$L$;
\item
$L(q)$, $q\in\Q$, is the abelian group~$L$ equipped with the symmetric
bilinear form
$x\otimes y\mapsto q(x\cdot y)$,
provided that it is still a lattice;
\item
$qL\subset L\otimes\Q$, $q\in\Q$, is the subgroup $\{qx\,|\,x\in L\}$,
also equipped with the restricted bilinear form;
as an abstract lattice, $qL\cong L(q^2)$.
\endroster
The same notation applies to discriminant forms whenever it makes sense.
Note that, if $p$ is a prime, $L(\frac1p)$ is a lattice if and only if
$\ell(\discr_pL)=\rank L$; this lattice is even if and only if
$L$ is even and either $p\ne2$ or
$\discr_2L$ is even.

Usually, we do \emph{not} assume isometries bijective;
for an isometry $\psi\:L\to S$, one has $\Ker\psi\subset\ker L$.
The group of bijective
auto\-isometries of a lattice~$L$ is denoted by $\OG(L)$.
There is a canonical homomorphism $\OG(L)\to\Aut\discr L$.

A \emph{$4$-polarization} of a lattice~$L$ is a distinguished vector $h\in L$ of
square~$4$; this vector is usually assumed but not present in the notation.
The group of polarized autoisometries is denoted by $\OG_h(L)$.
A \emph{line} in a $4$-polarized hyperbolic lattice~$L$ is an element of the
set
\[*
\Fn L:=\bigl\{a\in L\bigm|a^2=-2,\ a\cdot h=1\bigr\}.
\]
The set $\Fn L$ is finite; it admits a natural action of $\OG_h(L)$.

Two lattices $L'$, $L''$ are said to be in the same \emph{genus} if
$L'\otimes\R\cong L''\otimes\R$ and $L'\otimes\Q_p\cong L''\otimes\Q_p$ for
each prime~$p$. Each genus contains finitely many isomorphism classes.
According to~\cite{Nikulin:forms}, the genus of an \emph{even} nondegenerate
lattice is determined by its rank, signature, and discriminant form.
A realizability criterion is given by the following theorem.

\theorem[{see \cite[Theorem 1.10.1]{Nikulin:forms}}]\label{th.Nikulin}
A nondegenerate even lattice~$L$ with given inertia indices $(\sigma_+,\sigma_-)$
and discriminant form~$\CL$ exists if and only if
\roster
\item\label{Nikulin.length}
$\ell(\CL)\le r:=\rank L=\sigma_++\sigma_-$,
\item\label{vdB}
$\Br\CL=\sigma_+-\sigma_-\bmod8$ \rom(van der Blij formula~\cite{vanderBlij}\rom),
\endroster
and the following conditions are satisfied\rom:
\roster*
\item
$\ls|\CL|\det_p\CL=(-1)^{\sigma_-}\bmod(\Z_p^\times\!)^2$ for any prime $p>2$
for which $\ell_p(\CL)=r$\rom;
\item
either $\ell_2(\CL)<r$, or $\CL_2$ is odd, or
$\ls|\CL|\det_2\CL=\pm1\bmod(\Z_2^\times\!)^2$.
\endroster
\endtheorem

We fix the following notation for a few special lattices:
\roster*
\item
$\bH_n:=\bigoplus_{i=1}^n\Z\e_i$, $\e_i^2=-1$; once the basis is fixed,
we have a distinguished characteristic vector
$\be:=\e_1+\ldots+\e_n\in\bH_n$;
\item
$\bU:=\Z u_1+\Z u_2$, $u_1^2=u_2^2=0$, $u_1\cdot u_2=1$, is the
\emph{hyperbolic plane};
\item
$\bA_n$, $\bD_n$, $\bE_n$ are the \emph{negative} definite
lattices generated
by the root systems of the same name, see~\cite{Bourbaki:Lie};
\item
$\bL:=H_2(X)\cong\bE_8^2\oplus\bU^3$ is the intersection form of a
$K3$-surface~$X$ over~$\C$;
\item
$\bS_{p,\Artin}:=\NS(X)$
is the N\'{e}ron--Severi lattice of a supersingular $K3$-surface~$X$ over a field
of characteristic~$p$ with Artin invariant $\Gs=1,\ldots,10$.
\endroster
Recall that $\bA_n$ can be interpreted as the orthogonal complement
$\be^\perp\subset\bH_{n+1}$ and $\bD_n$ is the maximal even sublattice
in~$\bH_n$.
The nondegenerate even lattice $\bS_{p,\Artin}$ is
uniquely determined by
the properties $\Gs_+\bS_{p,\Artin}=1$, $\Gs_-\bS_{p,\Artin}=21$, and
$\discr\bS_{p,\Artin}$ is a $p$-elementary group of length $2\Artin$,
even if $p=2$. Similarly, $\bL$ is the only even unimodular
lattice with $\Gs_+\bL=3$ and $\Gs_-\bL=19$.

We also use freely the classification of definite unimodular lattices of
small rank found in~\cite{Conway.Sloane},
explaining the extra notation $L^+$ on the fly: usually, it stands for the
only ``interesting'' unimodular extension of~$L$.

\subsection{Lattice extensions\pdfstr{}{ {\rm(see~\cite{Nikulin:forms})}}}\label{s.extension}
From now on, unless specified otherwise, all lattices considered are even and
nondegenerate. Respectively, $q$-isotropic subgroups of a finite quadratic
form~$\CL$ are called just isotropic.

An \emph{extension} of a lattice~$S$ is any overlattice $L\supset S$. Two
extensions $L',L''\supset S$ are \emph{isomorphic} if there is a bijective
isometry $L'\to L''$ identical on~$S$.
One can also fix a subgroup $G\subset\OG(S)$ and speak about
\emph{$G$-isomorphisms} of extensions, \ie, bijective isometries whose
restriction to~$S$ is in~$G$.

Let $L\supset S$ be a finite index extension. Then we have natural inclusions
\[*
S\subset L\subset L\dual\subset S\dual,
\]
\cf.~\eqref{eq.dual}, and, hence, a well defined subgroup
$\CK:=L/S\subset\discr S=S\dual\!/S$. This subgroup is isotropic (since $L$
is an even integral lattice); it is called the \emph{kernel} of the extension
$L\supset S$. Conversely, if $\CK\subset\discr S$ is isotropic,
the lattice
\[*
L:=\bigl\{x\in S\otimes\Q\bigm|x\bmod L\in\CK\bigr\}
\]
is an extension of~$S$.
(If $\CK$ is $b$-, but not $q$-isotropic, then $L$ is odd.)
We say that $L$ is the extension of~$S$ by~$\CK$ (or by any collection of
vectors $a_1,a_2,\ldots\in S\otimes\Q$ such that $a_i\bmod S$ generate~$\CK$).
Thus, we have the following statement.

\theorem[{see \cite{Nikulin:forms}}]\label{th.extension}
Given a subgroup $G\subset\OG(S)$, the map
\[*
(L\supset S)\longmapsto\CK:=L/S\subset\discr S
\]
is a one-to-one
correspondence between the set of $G$-isomorphism classes of finite index
extensions $L\supset S$ and the set of
$G$-orbits of isotropic subgroups $\CK\subset\discr S$.
Under this correspondence, one has $\discr L=\CK^\perp\!/\CK$.
\endtheorem

In general, an extension $L\supset S$ can be described by fixing a finite
index sublattice $T\subset S_L^\perp$: then $L$ is a finite index extension
of $S\oplus T$ and, as such, is determined by an isotropic subgroup
\[*
\CK\subset\discr(S\oplus T)=\discr S\oplus\discr T.
\]
This subgroup~$\CK$ can be regarded as the graph of an anti-isometric additive
relation (also known as partially defined multi-valued homomorphism)
\[*
\psi\:\discr S\dashrightarrow\discr T;
\]
denoting by $\pr_S$, $\pr_T$
the projections to the two summands, we have
\[*
\alignedat2
\operatorname{Domain}\psi&=\pr_S(\CK),&\qquad
 \Ker\psi&=\CK\cap\discr S,\\
\Im\psi&=\pr_T(\CK),&\qquad
 \operatorname{Indet}\psi&=\CK\cap\discr T.
\endalignedat
\]
Hence, if $T=S_L^\perp$ is primitive in~$L$, then
$\psi\:\pr_S(\CK)\to\discr T$ is a conventional anti-isometry; if $S$ is also
primitive, then $\psi$ is injective.
With $T$ fixed and $G\subset\OG(S)$ as above, the $G$-isomorphism classes of
extensions are enumerated by the orbits of the two-sided action of
$G\times\OG(T)$ on the set of anti-isometric additive relations~$\psi$.
Note, though, that if we do not insist that $T$ should be primitive, distinct
pairs $(T,\psi)$ may give rise to isomorphic extensions.

An important consequence is the following restriction on the genus of~$T$.

\proposition[see \cite{Nikulin:forms}]\label{prop.genus.T}
If both~$S$ and~$T$ are primitive in~$L$ and $\discr_pL=0$ for some
prime~$p$, then $\psi_p\:\discr_pS\to\discr_pT$ is a bijective anti-isometry.
\endproposition

\subsection{Lemmas on discriminant forms}\label{s.discr.lemmas}
In this section, we state a few lemmas which would help us identify negative
definite lattices.

\lemma\label{lem.2-discr}
Let $T$ be a lattice with $2$-elementary group $\discr_2T$. Then
there is a finite index sublattice $T'\subset T$
such that $\discr_pT'=\discr_pT$ for all
primes $p\ne2$ and
$\discr_2T'$ is a $2$-elementary group of maximal length\rom:
$\ell(\discr_2T')=\rank T'$.
\endlemma

\proof
We start with the sublattice $2T$ and extend it to $T_0\subset T$ \via\ the
obviously
isotropic subgroup $4\discr_2(2T)$; the new discriminant $\CT:=\discr_2T_0$
has only $2$- and $4$-torsion. Such discriminant forms have been studied
in~\cite{DIK}. There is a well-defined nondegenerate
symmetric bilinear form
\[*
\cc\:2\CT\otimes2\CT\longto\Q/\Z,\quad
 2x\otimes2y\longmapsto2(x\cdot y),
\]
and, given an isotropic subgroup $\CK\subset\CT$,
one has (see \cite[4.2.2]{DIK})
\[*
2(\CK^\perp\!/\CK)=(\CK\cap2\CT)^\perp_{\cc}\big/
 \bigl(\CK\cap(\CK\cap2\CT)^\perp_\cc\bigl);
\]
in particular, $\CK^\perp\!/\CK$ is $2$-elementary if and only if
$(\CK\cap2\CT)^\perp_\cc$ is $\cc$-isotropic.
The kernel of the extension $T\supset T_0$ is $2\CT$;
hence, a
lattice~$T'$ as in the statement is obtained by extending~$T_0$ by any maximal
$\cc$-isotropic subgroup $\CK'\subset2\CT$.
Note that we always have the congruence
$\ell(2\CT)=\rank T-\ell(\discr_2T)=0\bmod2$ and any nondegenerate
symmetric bilinear form on an $\F_2$-vector space of even dimension is
null-cobordant (as follows from the classification of such forms).
\endproof

The arguments of~\cite{DIK} can easily be extended to finite forms with
$3$- and $9$-torsion only. Given such a form~$\CT$, we have a well-defined
nondegenerate quadratic form
\[*
\qc\:3\CT\longto\Q/2\Z,\quad 3x\longmapsto3x^2.
\]
The following statement is immediate.

\lemma\label{lem.3-forms}
Given a quadratic form~$\CT$ as above and an isotropic subgroup
$\CK\subset\CT$, one has $\ell(\CK^\perp\!/\CK)=\ell(\CT)$ if and only if
$\CK\subset3\CT$ and $\CK$ is $\qc$-isotropic.
\endlemma

Furthermore, for a $3$-elementary quadratic form~$\CT$, using the obvious
additivity, one can easily check the congruences
\[*
\Gd:=\Br\CT-2\ell(\CT)=0\bmod4,\qquad
\det\CT=(-1)^{\delta/4}\cdot\ls|\CT|\1.
\]
These congruences apply to~$\qc$; combining, for a quadratic
form~$\CT$
with $3$- and $9$-torsion only, we have
\[
\Gd:=\Br\CT+\Br\qc-2\ell(\CT)=0\bmod4,\qquad
\det\CT=(-1)^{\delta/4}\cdot\ls|\CT|\1.
\label{eq.3-forms}
\]
(Recall that $\Br q=0$ for any form~$q$ on $\CG9$.)

\lemma\label{lem.2-3-discr}
Let $T$ be a negative definite lattice with the following properties\rom:
\roster*
\item
$\discr_2T$ is odd $2$-elementary and $\Br(\discr_2T)=\rank T\bmod8$,
\item
$\discr_3T$ is $3$-elementary,
and
\item
$\discr_pT=0$ for all primes $p>3$.
\endroster
Then $T$ contains a
finite index sublattice $T'\cong\barT(6)$, where $\barT$ is odd
unimodular.
\endlemma

\proof
As in the proof of \autoref{lem.2-discr}, we extend $3T$ \via\
$9\discr_3(3T)$ to obtain a sublattice $T_0\subset T$ whose discriminant
$\discr_3T_0$ has only $3$- and $9$-torsion. Using
\autoref{th.Nikulin}\iref{vdB} and~\eqref{eq.3-forms}, we obtain
$\Br\qc=0$. Hence, $\qc$ is null-cobordant and
$T$ contains a sublattice~$T_3$ with $3$-elementary group
$\discr_3T_3$ of maximal length.

There remains to apply \autoref{lem.2-discr} to~$T_3$ to produce a sublattice
$T'\subset T$ with both $2$- and $3$-discriminants elementary and of maximal
length. Then $\barT:=T'(\frac16)$ is integral and unimodular;
it is odd since so is $\discr_2T'$.
\endproof

The next lemma describes a maximal $3$-elementary finite index sublattice.

\lemma\label{lem.3-discr}
Any
$3$-elementary lattice~$T$ has a finite index sublattice
$T'=\barT\dual\!(3)$, where $\discr\barT=\<\frac23r\>^m$ and
$r=\pm1$, $m\le2$ are such that $2mr=\Gs(T)\bmod8$.
\endlemma

\proof
Consider the extension $T_0\supset3T$ by $9\discr_3(3T)$,
followed by the extension $T_1\supset T_0$ by any maximal isotropic subgroup
of $(3\discr T_0,\qc)$.
Then, $T_1=\barT(3)$, where $\barT$ is as in the statement, and
$T'\supset T_1$ is the extension by the isotropic subgroup $3\discr T_1$.
\endproof

%
%

The following well-known lemma is easily proved by induction.

\lemma\label{lem.counts}
Let $\CA$ be an affine subspace in a quadratic $\F_3$-vector space, and let
$n_r(\CA)$ be the number of vector in~$\CA$ of square
$\frac23r$, $r\in\F_3$.
Then, for each $r\in\F_3$, one has either $\dim\CA\le2$ and $n_r(\CA)\le5$ or
$n_r(\CA)=0\bmod3$.
\endlemma

Finally, consider a finite quadratic form~$\CS_n$ generated by $n$
orthogonal elements $\Ga_i$, each of order~$2$ and square $\frac12$. We have
an obvious inclusion $\SG{n}\subset\Aut\CS_n$, the symmetric group acting
\via\ permutations of the generators.
The \emph{reflection} against an element $\Ga\in\CS_n$, $\Ga^2=1$, is
the autoisometry $t_\Ga\:x\mapsto x-2(x\cdot\Ga)\Ga$.
The group $\Aut\CS_n$ is generated by all reflections, whereas
$\SG{n}$
is generated by the reflections against vectors of Hamming norm~$2$.
Denote by $\rt_i\in\Aut\CS_n$ the reflection
against $\Ga_i+\ldots+\Ga_{i+5}$ (assuming that $n\ge i+5$), and
let $\rt\in\Aut\CS_n$ be the reflection against $\Ga_1+\ldots+\Ga_{10}$ (assuming that
$n\ge10$).

\lemma[D.~Pasechnik, private communication]\label{lem.Dima}
For $n\le9$, a complete list of representatives of the double cosets
$\frak{S}_n:=\SG{n}\!\backslash\!\Aut\CS_n/\SG{n}$ is as follows\rom:
\[*
\operatorname{identity},\quad
\rt_1\ \text{\rom(if $n\ge6$\rom)},\quad
\rt_1\rt_3\ \text{\rom(if $n\ge8$\rom)},\quad
\rt_1\rt_4\ \text{\rom(if $n\ge9$\rom)}.
\]
The set $\frak{S}_{10}$
is represented by $\{\rt_1\rt_3\rt_5\}\cup\{u,\rt u\,|\,u\in\frak{S}_9\}$\rom;
one has $\ls|\frak{S}_{10}|=9$.
\endlemma

In practice, we use \autoref{lem.Dima} to classify the bijective
anti-isometries between two copies, $\CS_n$ and~$-\CS_n$, each equipped with
a basis canonical up to order, up to the two-sided action of the group
$\SG{n}\times\SG{n}$, \cf. \autoref{s.extension}.
To do so, we identify the two groups by means of some bijection of their
bases and, with a certain abuse of the language, speak about the
anti-isometries~$\rt_i$, $\rt$, \etc.

\subsection{Generalized quadrangles\pdfstr{}{ {\rm(see~\cite{Payne.Thas})}}}\label{s.GQ}
The intersection of a quartic~$X$ and a plane in $\Cp3$ is a curve of
degree~$4$. It may happen that this curve is completely reducible, \ie,
splits into four lines $l_1,\ldots,l_4$.
(Note that these lines must be pairwise distinct, as otherwise $X$ would have
a singular point.)
If this is the case, we say that the
lines $l_1,\ldots,l_4$ constitute a \emph{plane} $\Ga\subset\Fn X$, \cf.
\autoref{lem.plane} below and the definition thereafter.
By definition, each plane consists of four lines. Occasionally, we consider
subconfigurations $\CF\subset\Fn X$ with the following properties:
\roster*
\item
each line is contained in a certain fixed number $p\ge2$ of planes;
\item
if two lines $l_1,l_2\in\CF$ intersect, they are contained in a plane
$\Ga\subset\CF$.
\endroster
In this case, renaming (lines, planes) to (points, lines) and taking
the inclusion
for the incidence relation, we obtain a combinatorial structure known as a
\emph{generalized quadrangle}
of order $(3,p-1)$, or $\GQ(3,p-1)$.
Specifically, a $\GQ(s,t)$ consists of two sets, $\Cal P$ (points) and
$\Cal B$ (lines), and an incidence relation~$|$, so that
\roster*
\item
each point is incident with $1+t\ge2$ lines,
\item
each line is incident with $1+s\ge2$ points, and
\item
for a point~$l$ and line $\Ga\notdivides|l$, there is a unique pair $\Ga',l'$
such that $l\divides|\Ga'\divides|l'\divides|\Ga$,
\endroster
see~\cite{Payne.Thas} for the
precise definition and further details.
(For the last axiom, one uses the
obvious fact that
a line $l\in\Fn X$ not contained in a plane $\Ga$ intersects exactly
one line $l'\in\Ga$, \cf. \autoref{rem.rels} below;
then, $l,l'$ are contained in a plane~$\Ga'$.)

According
to \cite[\S6.2]{Payne.Thas}, a generalized quadrangle $\GQ(3,t)$
exists if and only if $t=1,3,5,9$, and,
unless $t=3$, a quadrangle
is unique up to isomorphism. In the exceptional case $t=3$,
there are two quadrangles: $Q(4,3)$ and its dual $W(3)$.
Here, $Q(d,q)$, $d=3,4,5$,
can be described as the collection of points and lines in
a fixed nonsingular quadric of index~$2$ in the projective space $\Cp{d}$
over $\F_q$, whereas $W(q)$ is
the collection of points in the projective space $\Cp3$ over $\F_q$ and lines
Lagrangian with respect to any fixed symplectic form.

By our definition, two lines in a
generalized quadrangle $\CF\subset\Fn X$ intersect if and only if they are
contained in a plane $\Ga\subset\CF$. Hence, the adjacency graph of the lines
is uniquely determined by the combinatorics, and we denote by
$\bQ_{16}$, $\bQ_{40}'$, $\bQ_{40}''$, $\bQ_{64}$, $\bQ_{112}$
the corresponding lattices modulo kernel,
with the $4$-polarization defined as the sum of the
four lines constituting any plane (\cf. \autoref{lem.plane} below).
(By convention, the lattices $\bQ_{40}'$ and $\bQ_{40}''$ correspond to
the generalized quadrangles
$Q(4,3)$ and $W(3)$, respectively.)
All five lattices are hyperbolic; we have
\roster*
\item
$\rank\bQ_{112}=22$ and $\ls|\Fn\bQ_{112}|=112$
(see \autoref{rem.Fermat} and \autoref{ss.GQ9}),
\item
$\rank\bQ_{64}=19$ and $\ls|\Fn\bQ_{64}|=64$
(see \autoref{ss.GQ5}),
\item
$\rank\bQ_{40}^*=16$ and $\ls|\Fn\bQ_{40}^*|=40$
(see \autoref{proof.char=2} and \autoref{ss.GQ3}),
\item
$\rank\bQ_{16}=10$ and $\ls|\Fn\bQ_{16}|=16$
(see \autoref{ss.GQ1}).
\endroster
(References indicate parts of the paper where the realizability of the
generalized quadrangles by configurations of lines in nonsingular quartics is
discussed.)

\section{$K3$-surfaces}\label{S.K3}

Here, we give a brief account of the theory of $K3$-surfaces; for more
details and further references, we address the reader to~\cite{Huybrechts}.

\subsection{$K3$-surfaces}\label{s.K3}
An \emph{\rom(algebraic\rom) $K3$-surface} over an
algebraically closed field~$\Bbbk$ is a complete
nonsingular variety~$X$ over~$\Bbbk$ of dimension two such that
\[*
\Omega^2_X\cong\CO_X,\quad
H^1(X;\CO_X)=0.
\]
If $X$ is a $K3$-surface, the canonical epimorphism
$\operatorname{Pic}X\onto\NS(X)$ is an isomorphism;
furthermore, the lattice $\NS(X)$ is even and hyperbolic and
$\rank\NS(X)\le22$.

If $\Bbbk=\C$, one also considers \emph{analytic $K3$-surfaces}, which are
simply connected compact complex surfaces with the trivial canonical bundle.
All $K3$-surfaces are K\"{a}hler.
In general, $\Gs_+\NS(X)\le1$, and $X$ is algebraic if and only if
$\Gs_+\NS(X)=1$; in this case, $\NS(X)$ is nondegenerate.
We have a primitive embedding
\[*
\NS(X)\subset H_2(X;\Z)\cong\bL=\bE_8^2\oplus\bU^3,
\]
see \autoref{s.lattice};
hence, $\rank\NS(X)\le20$.
These statements on $\NS(X)$
extend to $K3$-surfaces over any
algebraically closed field~$\Bbbk$ of characteristic~$0$.

A $K3$-surface~$X$ is called \emph{\rom(Shioda\rom) supersingular} if
$\rank\NS(X)=22$.

\theorem[see \cite{Artin:supersingular}]\label{th.NS}
Assume that a $K3$-surface~$X$ over an algebraically closed field~$\Bbbk$ is
supersingular. Then $p:=\kchar>0$ and $\NS(X)\cong\bS_{p,\Artin}$ \rom(see
\autoref{s.lattice}\rom) for some integer $\Artin=1,\ldots,10$, called the
\emph{Artin invariant} of~$X$.
\endtheorem

If $X$ is not supersingular, then $\rank\NS(X)\le20$. Furthermore, according
to the next theorem, in this case we have (al least) the same restrictions on
$\NS(X)$ as in the case of characteristic~$0$.

\theorem[see~\cite{Deligne,Lieblich.Maulik}]\label{th.char=0}
If
a $K3$-surface~$X$ is not supersingular, there exists a
$K3$-surface~$X_0$ over a
field $\Bbbk_0$, $\fchar\Bbbk_0=0$, with the property that
$\NS(X_0)\cong\NS(X)$.
In particular, there exists a primitive extension $\bL\supset\NS(X)$.
\endtheorem

\subsection{Quartics}\label{s.quartics}
Any nonsingular quartic $X\subset\Cp3$ is a $K3$-surface. This surface is
equipped with a canonical $4$-polarization $h\in\NS(X)$, \viz. the hyperplane
sections; this polarization is always assumed when we speak about quartics.

\theorem[see \cite{Saint-Donat}]\label{th.Saint-Donat}
The
$4$-polarization $h\in\NS(X)$ of a nonsingular quartic $X\in\Cp3$
has the following property\rom:
there is no vector $e\in\NS(X)$ such that either
\roster
\item\label{Saint-Donat.exceptional}
$e^2=-2$ and $e\cdot h=0$ \rom(\emph{exceptional divisor}\rom) or
\item\label{Saint-Donat.double}
$e^2=0$ and $e\cdot h=2$.
\endroster
Conversely, given a $K3$-surface~$X$
over
an algebraically closed field~$\Bbbk$, $\fchar\Bbbk\ne2$,
a
$4$-polarization $h\in\NS(X)$ contained in the positive cone
of $\NS(X)$
and
satisfying
the two conditions
above
embeds~$X$ into~$\Cp3$ as a nonsingular
quartic.
\endtheorem

\remark\label{rem.Saint-Donat}
Strictly speaking,
\cite{Saint-Donat} makes a global assumption that $\fchar\Bbbk\ne2$.
Nevertheless, the arguments leading to the necessity of
conditions~\iref{Saint-Donat.exceptional}, \iref{Saint-Donat.double} work in
any characteristic, as essentially they rely upon the Riemann--Roch theorem
only. It would not be a surprise if the conditions were also sufficient.
Unfortunately, I could not find in the literature an appropriate reference
covering all the details of the proof;
that is why I
refrain from stating the realizability of the found configurations in
\autoref{th.char=2}, leaving the latter as an ``upper bound.''
\endremark

A $4$-polarized hyperbolic lattice~$S$ satisfying the necessary
conditions~\iref{Saint-Donat.exceptional}, \iref{Saint-Donat.double}
in \autoref{th.Saint-Donat} is called \emph{admissible}.

A \emph{geometric realization} of an admissible
lattice~$S$ is a lattice
extension $\bL\supset S$
or $\bS_{p,\Artin}\supset S$ (see \autoref{s.lattice}),
where $1\le\Artin\le10$ and $p$ is a prime;
we also require that
\roster*
\item
the primitive hull $\tS:=(S\otimes\Q)\cap\bL$ (in the former case) or
\item
the $4$-polarized lattice $(\bS_{p,\Artin},h)$ (in the latter case)
\endroster
should
still
be admissible.
In view of Theorems~\ref{th.NS} and~\ref{th.char=0},
the following simple consequence of \autoref{th.Nikulin}
(applied to the orthogonal complement of $T:=S^\perp$ in~$\bL$
or~$\bS_{p,\Artin}$)
and
\autoref{prop.genus.T} gives us a necessary condition for the existence of a
primitive geometric realization of a given hyperbolic lattice~$S$.

\theorem\label{th.realization}
Consider a primitive
hyperbolic sublattice
$S\subset\NS(X)$ and denote $\delta:=22-\rank S$ and $\CS:=\discr S$.
If $X$ is supersingular, let $p:=\kchar$\rom; otherwise, let $p:=0$.
Then we have $\ell_q(\CS)\le\delta$ for each prime $q\ne p$ and
\roster*
\item
$\ls|\CS|\det_q\CS=1\bmod(\Z_q^\times\!)^2$ for any prime
$q\ne2$ or~$p$ for which $\ell_q(\CS)=\delta$\rom;
\item
either $p=2$, or $\ell_2(\CS)<\delta$, or $\CS_2$ is odd, or
$\ls|\CS|\det_2\CS=\pm1\bmod(\Z_2^\times\!)^2$.
\endroster
\endtheorem

(Note that,
by the additivity of both Brown invariant~$\Br$ and
signature $\Gs_+-\Gs_-$,
condition~\iref{vdB} in \autoref{th.Nikulin} holds automatically.)
An important observation is the fact that the restriction imposed by
\autoref{th.realization} at a prime $q\ne\kchar$ does not depend on $\kchar$.

If $\kchar=0$, the condition given by \autoref{th.realization} is
also sufficient, as follows from the surjectivity of the period
map~\cite{Kulikov:periods}
(see also~\cite{DIS} for the details on the moduli space).
If $X$ is supersingular, a sufficient condition
is that at least one of the (finite set of)
extensions $\bS_{p,\Artin}\supset S$ obtained by
Nikulin's construction (see \autoref{s.extension}) should be admissible.
Note, though, that in this case $\Fn S$ may be a proper subset of the set
$\Fn X=\Fn\bS_{p,\Artin}$; these phenomena are discussed in \autoref{S.exotic}
and \autoref{S.2-exotic}.

Let $X\subset\Cp3$ be a nonsingular quartic. Then, sending a line $l\subset X$ to
its class $[l]\in\NS(X)$, we obtain a map $\Fn X\to\NS(X)$.

\lemma[\cf. \cite{DIS}]\label{lem.lines}
The map $l\mapsto[l]$ establishes a
bijection
$\Fn X=\Fn\NS(X)$.
\endlemma

\proof
The map is obviously well defined and injective (since each line
has negative self-intersection), and its image is in $\Fn\NS(X)$.
By the Riemann--Roch theorem, any element
$a\in\Fn\NS(X)$ is realized by a unique $(-2)$ curve~$C$.
Assume that $C$ is reducible, $C=C_1+\ldots+C_k$, where all components~$C_i$
are also $(-2)$-curves.
Then, since $1=C\cdot h=\sum_iC_i\cdot h$,
all but one
components of~$C$ are exceptional divisors, contradicting
\autoref{th.Saint-Donat}\iref{Saint-Donat.exceptional}.
Thus, $C$ is irreducible;
since also $C$ has projective degree $1=C\cdot h$,
we conclude that $C$ is a line.
\endproof

\subsection{(Quasi-)elliptic pencils}\label{s.elliptic}
Let
$\pi\:X\to\Cp1$ be a pencil of curves of
arithmetic genus~$1$. If a generic fiber of~$\pi$ is a smooth elliptic curve,
the pencil~$\pi$ is \emph{elliptic}; otherwise
(generic fiber is singular), $\pi$ is
\emph{quasi-elliptic}.

All fibers of a (quasi-)elliptic pencil $\pi\:X\to\Cp1$ are linearly
equivalent and, for each fiber~$F$, one has $F^2=0$. Conversely,
if $X$ is a $K3$-surface, then each primitive,
effective, and numerically effective divisor $F\subset X$ such that $F^2=0$ is a
fiber of a unique (quasi-)elliptic pencil $\pi\:X\to\Cp1$.

Let $F=\sum_ir_iC_i$ be a reducible fiber. Each reduced component~$C_i$ is a
smooth rational $(-2)$-curve, and the dual intersection graph of~$F$ is a
certain affine Dynkin diagram~$\aD$, see~\cite{Bourbaki:Lie}.
Denoting by $\Z\aD$ the
intersection lattice freely generated by the vertices $C_i\in\aD$,
the kernel $\ker\Z\aD$ has rank~$1$ and is generated by a unique positive linear
combination $\sum_ir_iC_i$; this generator of $\ker\Z\aD$ is~$F$.

We define the \emph{Milnor number} of an (affine or elliptic) Dynkin
diagram~$D$ as the rank $\mu(D):=\rank(\Z D/\ker)$.
Thus, $\mu(D)$ is the number of vertices of~$D$ in the elliptic case and the
number of vertices minus~$1$ in the affine case.

\theorem[see \cite{Rudakov.Shafarevich}]\label{th.quasi}
Let $\pi\:X\to\Cp2$ be a pencil of curves of arithmetic genus~$1$,
and denote by $\aD_1,\ldots$ the
dual intersection graphs of the components of
the reducible fibers of~$\pi$.
Then, $\pi$ is quasi-elliptic if and
only if
\roster
\item\label{quasi.1}
$p:=\kchar=2$ or~$3$\rom;
\item\label{quasi.2}
each lattice $\Z\aD_i/\ker$ is $p$-elementary\rom;
\item\label{quasi.3}
one has $\sum_i\mu(\aD_i)=b_2(X)-2$.
\endroster
\endtheorem

If $X$ is a $K3$-surface, then $b_2(X)=22$ and $e(X)=12\chi(X)=24$ (the
\'{e}tale Euler characteristic).
Recall also that the affine Dynkin diagrams~$\aD$ with $p$-elementary lattice
$\Z\aD/\ker$ are:
\roster*
\item
\leavevmode\smash{$\tA_1$}, \smash{$\tE_7$}, \smash{$\tE_8$}, \smash{$\tD_{2k}$},
if $p=2$, and
\item
\leavevmode\smash{$\tA_2$}, \smash{$\tE_6$}, \smash{$\tE_8$}, if $p=3$.
\endroster
If a pencil $\pi\:X\to\Cp2$ is elliptic, instead of
\autoref{th.quasi}\iref{quasi.3} we have the identity
\[
\sum\bigl(e(F_i)+d(F_i)\bigr)=e(X),
\label{eq.Euler}
\]
the summation running over all singular fibers~$F_i$ of~$\pi$.
Here, $d(F_i)\ge0$ is the wild ramification index and
\[*
e(F_i)\ge\ls|\text{irreducible components of $F_i$}|\ge1.
\]

Let $X\subset\Cp3$ be a nonsingular quartic, and let $\pi\:X\to\Cp1$ be
a \hbox{(quasi-)}\allowbreak elliptic pencil.
A fiber of~$\pi$ consisting entirely of lines is called \emph{parabolic};
any other singular fiber is called \emph{elliptic}.
Each parabolic fiber is an affine Dynkin diagram $\aD\subset\Fn X$,
whereas the
configuration of lines contained in an elliptic fiber is a Dynkin diagram,
possibly empty or disconnected, $D\subset\Fn X$;
the type of~$D$ is referred to as the \emph{linear type} of the elliptic
fiber.
Regarded as spatial curves, all fibers of~$\pi$ have the same
degree. This fact limits the types of parabolic fibers and linear types
of
elliptic fibers appearing in the same pencil.

We denote by $\flines\pi$ the number of lines contained in the fibers
of~$\pi$. The following statement is an immediate consequence of
\autoref{th.quasi}\iref{quasi.3} and~\eqref{eq.Euler}.

\corollary\label{cor.count}
If a pencil $\pi\:X\to\Cp1$ is elliptic, then
\[*
\flines\pi\le\ls|\text{\rm components in the singular fibers of~$\pi$}|\le24;
\]
If $\pi$ is quasi-elliptic, then
\[*
\flines\pi\le20+\ls|\text{\rm parabolic fibers of $\pi$}|.
\]
\endcorollary

When applying the first inequality, we often use the fact that the upper
bound~$24$ is reduced by at least~$1$ by each elliptic fiber of~$\pi$, as
such a fiber contains at least one component that is not a line.

\section{Configurations and pencils}\label{S.config}

In this section, we discuss simplest arithmetical properties of configurations
of lines. Most results here either are contained in~\cite{DIS} or can be
regarded as immediate extensions/generalizations thereof.

\subsection{Configurations of lines}\label{s.config}

Let~$S$ be a $4$-polarized hyperbolic lattice; we will always assume
that $S$ is admissible, \ie, satisfies
the conditions in \autoref{th.Saint-Donat}.
Recall (see \autoref{s.lattice})
that the \emph{configuration of lines} in~$S$ is the set
\[*
\Fn S:=\bigl\{a\in S\bigm|a^2=-2,\ a\cdot h=1\bigr\}.
\]
Elements of~$\Fn S$ are called \emph{lines}.
The hyperbolicity of~$S$ and \autoref{th.Saint-Donat} imply that, for any
pair $a_1,a_2\in\Fn S$, one has $a_1\cdot a_2=0$ or~$1$
(\cf.~\cite{DIS}); respectively, we say
that the two lines $a_1,a_2$ are \emph{disjoint} or \emph{intersect}.
The set $\Fn X$ is often treated as a graph, with two lines
connected by an edge if and only if they intersect.

The next few lemmas are based on the following simple observation,
\cf.~\cite{DIS}: each sublattice $S'\subset S$ containing~$h$ and at least
one line is hyperbolic; hence, $\ker S'=0$ and any numeric relation
$u=0\bmod(S')\dual$ implies a true relation $u=0$.

\lemma\label{lem.plane}
Assume that four lines $a_1,\ldots,a_4\in\Fn S$ intersect each other, \ie,
\[*
a_i\cdot a_j=1\ \text{for $1\le i<j\le4$}.
\]
\rom(In other words, the lines constitute the complete graph $K_4$.\rom)
Then
\[*
a_1+a_2+a_3+a_4=h.
\]
If $S=\NS(X)$, then the lines $a_1,\ldots,a_4$ are cut on~$X$ by a
plane.
\endlemma

A quadruple $a_1,\ldots,a_4\in\Fn S$ as in \autoref{lem.plane} is called a
\emph{plane} in~$S$.

The \emph{valency} $\val l$ of a line $l\in\Fn X$ is the number of lines
$a\in\Fn X$ intersecting~$l$
(alternatively, this is the valency of~$l$ as a vertex of the graph $\Fn X$),
whereas the \emph{multiplicity} $\mult l$ is
the number of planes $\Ga\subset\Fn X$ containing~$l$.

\lemma\label{lem.4-4}
Assume that eight lines $a_1,\ldots,a_4,b_1,\ldots,b_4\in\Fn S$
\figure
\cpic{4=4}
\caption{The configuration in \autoref{lem.4-4} (the graph $K_{4,4}$)}\label{fig.4-4}
\endfigure
intersect as shown in \autoref{fig.4-4}, \ie,
\[*
a_i\cdot a_j=b_i\cdot b_j=0\ \text{for $1\le i<j\le4$},\quad
a_i\cdot b_j=1\ \text{for all $i,j=1,\ldots,4$}.
\]
\rom(In other words, the lines constitute the complete bipartite graph
$K_{4,4}$.\rom) Then
\[*
a_1+a_2+a_3+a_4+b_1+b_2+b_3+b_4=2h.
\]
If $S=\NS(X)$, then the lines $a_1,\ldots,b_4$ are cut on~$X$ by a
quadric.
\endlemma

\lemma\label{lem.10-2}
Assume that twelve lines $a_1,\ldots,a_{10},b_1,b_2\in\Fn S$
\figure
\cpic{10=2}
\caption{The configuration in \autoref{lem.10-2} (the graph $K_{10,2}$)}\label{fig.10-2}
\endfigure
intersect as shown in \autoref{fig.10-2}, \ie,
\[*
a_i\cdot a_j=b_1\cdot b_2=0\ \text{for $1\le i<j\le10$},\quad
a_i\cdot b_j=1\ \text{for $i=1,\ldots,10$, $j=1,2$}.
\]
\rom(In other words, the lines constitute the complete bipartite graph
$K_{10,2}$.\rom) Then
\[*
a_1+a_2+\ldots+a_{10}+3b_1+3b_2=4h.
\]
\endlemma

In \autoref{prop.10-2.char2} below we show
that a configuration as in \autoref{lem.10-2} cannot exist if
$\kchar=2$.

\lemma\label{lem.6-2}
Assume that fourteen lines $a,b,a_1,\ldots,a_6,b_1,\ldots,b_6\in\Fn S$
\figure
\cpic{6=2}
\caption{The configuration in \autoref{lem.6-2}}\label{fig.6-2}
\endfigure
intersect as shown in \autoref{fig.6-2}, \ie,
\[*
a\cdot b=a\cdot b_i=a_i\cdot b=a_i\cdot a_j=b_i\cdot b_j=0,\
a\cdot a_i=b\cdot b_i=a_i\cdot b_i=1\
\]
for all pairs $i,j=1,\ldots,6$ such that $i\ne j$.
Then
\[*
3a+a_1+\ldots+a_6=3b+b_1\ldots+b_6.
\]
\endlemma

\remark\label{rem.rels}
The relations in Lemmas~\ref{lem.plane}--\ref{lem.6-2} can be used to assert
the existence and uniqueness of lines with certain properties. If all
but one lines of a configuration as in one of the lemmas are known,
there is at most one ``missing'' line completing the configuration; if
the corresponding coefficient in the relation equals $\pm1$, such a line does
exist (as it is a linear combination of the known ones).
In addition, we can use the lemmas to control the intersection of a line~$l$
other than the given ones with the lines constituting the configuration.
Thus,
\roster*
\item
in \autoref{lem.plane}, $l$ intersects exactly one of $a_1,\ldots,a_4$, and
\item
in \autoref{lem.4-4}, $l$ intersects exactly two of $a_1,\ldots,b_4$.
\endroster
(In the other two lemmas, one should weight the intersections with the
coefficients present in the relations.)
\endremark

\subsection{Pencils}\label{s.pencils}
The pencil of planes in~$\Cp3$ passing through a fixed line $l\subset X$
defines a (quasi-)elliptic pencil $\pi:=\pi[l]\:X\to\Cp1$:
the fibers of~$\pi$ are the
\emph{residual cubics} obtained by removing the
common component~$l$ from the
quartic curve cut on~$X$ by a plane. Clearly, the lines that are in the
fibers of~$\pi$ are precisely those intersecting~$l$. From this point of
view, there are two kinds of fibers: \emph{$3$-fibers}, split into three
lines (constituting, together with~$l$, a plane in $\Fn X$), and
\emph{$1$-fibers}, consisting of a single line (and an irreducible residual
conic). Each line in~$X$ that is disjoint from~$l$ intersects each fiber
of~$\pi$ at a single point; hence, it is a section of~$\pi$.

Motivated by this construction, given a line $l\in\Fn S$, we define the
\emph{maximal pencil} with the \emph{axis}~$l$ as the set
\[*
\pencil(l):=\bigl\{a\in\Fn S\bigm|a\cdot l=1\bigr\}.
\]
By \autoref{lem.plane}, the lines constituting $\pencil(l)$ split into
pairwise disjoint groups, each consisting of three or one line; they are
called $3$- and $1$-fibers of $\pencil(l)$, respectively.
The \emph{type} of a pencil is the pair $(p,q)$ representing the numbers of
its $3$- and $1$-fibers.
The \emph{pencil structure} of a configuration $\Fn S$ is the list of types
of all pencils $\pencil(l)$, $l\in\Fn S$; it is usually represented in the
partition notation and can be used as an easily comparable combinatorial
invariant.

Sometimes, we consider pencils $\pencil\subset\pencil(l)$ that are not
necessarily maximal; however, we always insist that, whenever $\pencil$
contains two lines $a_1$, $a_2$ that intersect, it also contains the third line
$h-l-a_1-a_2$ of the same $3$-fiber. In other words, $\pencil$ is a maximal
pencil in the sublattice $S'\subset S$ spanned by~$h$, $l$, and \emph{some}
of the lines $a\in\Fn S$ that intersect~$l$.

A \emph{section} of a pencil $\pencil\subset\pencil(l)$ is any line
$a\in\Fn S$ disjoint from~$l$. The number of sections of~$\pencil$ is denoted by
$\sect\pencil$.

Two pencils $\pencil(l_1)$, $\pencil(l_2)$ are called \emph{obverse} if
their axes $l_1$, $l_2$ are disjoint.

\subsection{The lattice $P\sb{p,q}$\pdfstr{}{ {\rm(see~\cite{DIS})}}}\label{s.P(p,q)}
Fix a pencil $\pencil\subset\Fn S$ (not necessarily maximal)
of type $(p,q)$ and denote by $P_{p,q}\subset S$ the sublattice spanned by
the polarization~$h$, axis~$l$, and the members of~$\pencil$.

Let
$\fb_3\pencil:=\{1,\ldots,p\}$ and $\fb_1\pencil:=\{1,\ldots,q\}$ be the
sets of the $3$- and $1$-fibers of~$\pencil$, respectively.
(We emphasize that we regard $\fb_3$ and $\fb_1$ as two \emph{disjoint} sets.)
The lattice $P_{p,q}$ is generated by the classes $h$, $l$,
$m_{i,\pm}$, $i\in\fb_3\pencil$ (two lines from each
$3$-fiber), and
$n_k$, $k\in\fb_1\pencil$ (the lines constituting the $1$-fibers).
The third line in the $i$-th $3$-fiber is $m_{i,0}=h-l-(m_{i,+}+m_{i,-})$, see
\autoref{lem.plane}. When speaking about dual vectors $h^*$, $l^*$, \etc., we
always refer to this distinguished basis (which involves some choice for each
$3$-fiber).

\observation\label{obs.discr}
The $3$-primary part $\discr_3P_{p,q}$ contains the classes
represented by the following mutually orthogonal vectors:
\roster*
\item
$\Gl:=\frac13(l-h)$: one has $\Gl^2=0$ and $\Gl\cdot h=\Gl\cdot l=-1$;
\item
$\Gm_i:=\frac13(m_{i,+}-m_{i,-})$, $i\in\fb_3\pencil$:
one has $\Gm_i^2=-\frac23$ and $\Gm_i\cdot h=0$.
\endroster
If $r:=p+q-1\ne0\bmod3$, then $\discr_3P_{p,q}$ is generated by $\Gm_i$,
$i\in\fb_3\pencil$, and
the order~$9$ class of the vector
\roster*
\item
$\Gu:=\frac13\bigl(l-r\Gl+\sum_{i=1}^p(m_{i,+}+m_{i,-})-\sum_{k=1}^qn_k\bigr)$;
\endroster
note that $3\Gu=-r\Gl\ne0\bmod P_{p,q}$.
Hence, in this case, the subgroup of elements of order~$3$ is
generated by~$\Gl$ and~$\Gm_i$.
If $p+q=1\bmod3$, then $\discr_3P_{p,q}$ is generated by $\Gl$, $\Gm_i$,
and the order~$3$ class of
\roster*
\item
$\Go:=\frac13\bigl(l+\sum_{i=1}^p(m_{i,+}+m_{i,-})-\sum_{k=1}^qn_k\bigr)$.
\endroster
The $\OG_h(P_{p,q})$-orbits of order~$3$ elements in $\discr_3P_{p,q}$ are
$r\oL:=\{r\Gl\}$, $r=\pm$, and
\[*
\textstyle
\oM_k:=\bigl\{r\Gl+\sum_{i\in I}\pm\Gm_i\bigm|
 r\in\F_3,\ I\subset\fb_3\pencil,\ \ls|I|=k\bigr\},\quad k=0,\ldots,p;
\]
if $p+q=1\bmod3$, there also are at most six orbits
\[*
\textstyle
r\oO_{s}:=\bigl\{\Ga\in\discr_3P_{p,q}\bigm|
 \Ga\cdot\Gl=-\frac13r,\ \Ga^2=\frac23s\bigr\},\quad r=\pm,\ s\in\F_3.
\]
Note also that any element of
$\Ker\bigl[\Aut\Gl^\perp\to\Aut(\Gl^\perp\!/\Gl)\bigr]$
lifts to $\OG_h(P_{p,q})$.

The $2$-primary part $\discr_2P_{p,q}$
is generated by the classes of $3\Gn_k$, where
\roster*
\item
$\Gn_k:=n_k^*=-\frac12(\Gl+n_k)$, $k\in\fb_1\pencil$:
one has $\Gn_k^2=-\frac12$ and $\Gn_k\cdot h=0$.
\endroster
With respect to this basis, the image of $\OG_h(P_{p,q})$
in $\Aut\discr_2P_{p,q}$ is the subgroup $\SG{q}$ permuting the generators,
and the orbits are
\[*
\textstyle
\oN_k:=\bigl\{\sum_{i\in I}\Gn_i\bigm|
 I\subset\fb_1\pencil,\ \ls|I|=k\bigr\},\quad k=0,\ldots,q.
\]
Note also that $\OG_h(P_{p,q})$ acts on the $2$- and $3$-torsion
independently, \ie, the image of $\OG_h(P_{p,q})$ in
$\Aut\discr P_{p,q}=\Aut\discr_2P_{p,q}\times\Aut\discr_3P_{p,q}$
is the product of its images in the two factors.
\endobservation

We are interested in the possible finite index extensions
$P_{p,q}\subset\tP\subset\NS(X)$; recall that each such extension is described
by its kernel $\tP/P_{p,q}$, which is an isotropic subgroup of
$\discr P_{p,q}$, see \autoref{th.extension}.
A pencil~$\pencil$ is called \emph{primitive} (in a given lattice~$S$) if
$\tP/P_{p,q}=0$; otherwise, it is called \emph{imprimitive}.
Due to \autoref{obs.discr}, the isotropic vectors are those in the orbits
$\pm\oL$, $\pm\oO_0$, $\oM_{3i}$, $i>0$, and $\oN_{4k}$, $k>0$.

\lemma\label{lem.discr}
The kernel $\tP/P_{p.q}$ is disjoint from
the orbits $\pm\oL$, $\oM_3$, and $\oN_4$.
If $p+q=1$, $4$, or~$7$, it is also disjoint from $\pm\oO_0$.
\endlemma

\proof
If $\Gl\in\NS(X)$, then so is $e:=-2\Gl$, which contradicts
\autoref{th.Saint-Donat}\iref{Saint-Donat.double}.

The orbits $\oM_3$ and $\oN_4$ are represented by the exceptional divisors
$\Gm_1+\Gm_2+\Gm_3$
and $\frac12(n_1-n_2+n_3-n_4)$, respectively,
see \autoref{th.Saint-Donat}\iref{Saint-Donat.exceptional}.
For the last statement, it suffices to consider the case $p=0$; then,
the vectors $2l^*$, $l^*+\Gn_1-\Gn_2-\Gn_3-\Gn_4$,
and $l^*-\sum_{k=1}^7\Gn_k$, respectively, are exceptional divisors.
\endproof

\subsection{Types of pencils}\label{s.types}
In the next statement, we show that, with one exception, a pencil of the form
$\pi[l]$, see \autoref{s.pencils}, is elliptic.

\proposition\label{prop.Euler}
With one exception, \viz. the case where
\roster*
\item
$\kchar=3$ and $\pencil(l)$ is of type $(10,0)$,
hence $\pi[l]$ is quasi-elliptic,
\endroster
the type $(p,q)$ of a pencil $\pencil\subset\NS(X)$
satisfies Euler's inequality
$3p+2q\le24$.
\endproposition

\proof
The fibers of $\pi[l]$ are all of type
$\tA_2$, $\tA_2^*$ ($p$ copies)
or $\tA_1$, $\tA_1^*$ ($q$ copies). By \autoref{th.quasi}, such a pencil can
be quasi-elliptic only if either
\roster*
\item
$\kchar=2$ and $(p,q)=(0,20)$, or
\item
$\kchar=3$ and $(p,q)=(10,0)$.
\endroster
The former possibility is eliminated by \autoref{th.realization} and
\autoref{lem.discr}, as the only isotropic orbits in $\discr_3P_{0,20}$
are $\pm\oL$, see \autoref{obs.discr}.
If $\pi[l]$ is elliptic, the inequality $3p+2q\le24$ follows
from \autoref{cor.count}.
\endproof

\proposition[see~\cite{DIS}]\label{prop.types}
Assume that $\kchar\ne2,3$ or $X$ is not supersingular. Then the type~$(p,q)$
of a pencil $\pencil(l)\in\Fn X$ takes values listed in \autoref{tab.Euler}.
\table
\caption{The bounds on $(p,q)$ and $\val l$, see \autoref{prop.Euler}}\label{tab.Euler}
\centerline{\setbox0\hbox{00}\vbox{\halign{
 \hss\strut$#$&&\quad\hbox to\wd0{\hss$#$\hss}\cr
p=&6&5&4&3&2&1&0\cr
q\le&2&3&6&7&9&10&12\cr
\val l\le&20&18&18&16&15&13&12\cr}}}
\endtable
Besides, if $p=6$ and $q>0$ or $(p,q)=(4,6)$, the pencil is necessary
imprimitive.
\endproposition

For a pencil~$\pencil$ of type~$(4,6)$, the imprimitivity implies
the existence of a section intersecting all ten fibers: by
\autoref{lem.discr}, the kernel $\tP/P_{4,6}$ is generated by~$\Go$, and the
section is given by \autoref{lem.10-2}. Pencils in supersingular quartics
over fields of characteristic~$2$ or~$3$ are considered in subsequent
sections; the restrictions are listed in \autoref{prop.2-exotic.counts} and
\autoref{cor.exotic.type}, respectively.

\proof[Proof of \autoref{prop.types}]
Since $3p+2q\le24$, we only need to eliminate the values $p=7$ and $(p,q)=(6,3)$
or $(5,4)$. To do so, assume that $P_{p,q}\subset\tP\subset\NS(X)$, where
$\tP$ is a finite index extension, and use \autoref{lem.discr} to bound the
size of the
group $\discr_3\tP$ from below. Then, \autoref{th.realization} implies that
$\tP$ does not admit a primitive geometric realization.
The imprimitivity is proved similarly, using the group
$\discr_3P_{p,q}$ instead of $\discr_3\tP$. For more details, see~\cite{DIS}.
\endproof

\proposition[see~\cite{DIS}]\label{prop.types.2}
Assume
that $\kchar\ne2$ or $X$ is not supersingular.
Then any pencil of type $(p,q)$, $p+q\ge11$, is imprimitive\rom; up
to reordering the $1$-fibers, the kernel
$\tP/P_{p,q}\subset\discr_2P_{p,q}$ is generated by
\roster*
\item
$3(\Gn_1+\ldots+\Gn_8)$ and $3(\Gn_5+\ldots+\Gn_{12})$ for $(p,q)=(0,12)$, or
\item
$3(\Gn_1+\ldots+\Gn_8)$ for $p+q=11$.
\endroster
\endproposition

This statement is proved similar to \autoref{prop.types}. As a consequence,
in the case $p+q=11$, each section intersects an even number of $1$-fibers
$n_1,\ldots,n_8$. If $q=12$, the fibers split into three groups,
$n_s,\ldots,n_{s+3}$, $s=1,5,9$, and the intersections of each section with
these groups are either all even or all odd.

Formally, \autoref{prop.types.2} still holds if $\kchar=2$ and $X$ is
supersingular, but the statement becomes void as $p+q\le8$ in this case (see
\autoref{prop.2-exotic} below).

\corollary\label{cor.types.2}
Assume that $\kchar\ne2$ or $X$ is not supersingular, and let $\pencil$ be a
pencil of type $(0,q)$. If $q\ge11$, a section~$s$ of~$\pencil$ can
intersect at most seven fibers of~$\pencil$\rom; if $q=12$, a section can
intersect at most six fibers.
\endcorollary

\proof
A section cannot intersect more than eight fibers by \autoref{lem.10-2}. If
$q=11$ and $s$ intersects eight fibers, then, up to reordering, these fibers
are either $n_1,\ldots,n_8$ or $n_1,\ldots,n_6,n_9,n_{10}$, see
\autoref{prop.types.2}. In the former case, $\tP$ contains a vector
as in \autoref{th.Saint-Donat}\iref{Saint-Donat.double}, \viz.
\[*
e:=-h+l+\tfrac12(n_1+\ldots+n_8)+s;
\]
in the latter case, it contains an
exceptional divisor $e-n_7-n_8$.
Similarly, if $q=12$ and $s$ intersects seven fibers, these fibers can be
chosen $n_1,n_2,n_3,n_5,n_6,n_7,n_9$; then, the vector
\[*
-h+l+\tfrac12(n_1+n_2+n_3-n_4+n_5+n_6+n_7-n_8)+s
\]
is an exceptional divisor.
\endproof

\section{Triangle free configurations}\label{S.trig.free}

A configuration $S$ is said to be \emph{triangle free} if the graph $\Fn S$
has no cycles of length~$3$. According to \autoref{lem.plane}, this
condition is equivalent to
the requirement that $S$ should contain no plane.

\subsection{Statement and setup}\label{s.trig.free}
The principal result of this section is the following theorem.
The proof follows that found in~\cite{DIS}. In
fact, in the present paper we are mainly interested in a few intermediate lemmas.

\theorem[see \autoref{proof.triangle.free}]\label{th.triangle.free}
If the configuration $\Fn X$ is triangle free and $\kchar\ne2$, then
$\ls|\Fn X|\le61$.
\endtheorem

\remark\label{rem.triangle.free}
Probably, the bound given by \autoref{th.triangle.free} is not sharp; the
best known example is a supersingular quartic~$X$ in characteristic~$7$
with $47$ lines:
\roster*
\item
$\ps=\PS{[ [ 0, 12 ], 20 ], [ [ 0, 10 ], 24 ], [ [ 0, 8 ], 3 ]}$,
$\Artin(X)=1$, and $\rank\Fano(X)=21$.
\endroster
The best (known to me) examples over other fields are $33$ lines if
$\Bbbk=\C$ and $37$ lines if $\kchar=3$ (see \autoref{prop.(0,q)} below).
\endremark

Throughout this section, we use the following simple observation. Let
$\aD\subset\Fn X$ be an affine Dynkin diagram, and let $\sum_ir_il_i$,
$l_i\in\aD$, be the positive generator of $\ker\Z\aD$. This generator,
regarded as a divisor, is obviously primitive,
effective, and numerically effective; hence,
the lines in~$\aD$ constitute
a reducible fiber of a \hbox{(quasi-)}\allowbreak elliptic pencil
$\pi:=\pi_{\aD}\:X\to\Cp1$, see \autoref{s.elliptic}.
Any other line $l\in\Fn X$ either is in a reducible singular
fiber of~$\pi$ or intersects each fiber of~$\pi$. Thus,
\[
\ls|\Fn X|=\flines\pi+
 \ls|\text{lines adjacent to a vertex of~$\aD$}|.
\label{eq.pencil.count}
\]

Another fact used freely without further references is that, in a triangle free
configuration $\Fn X$, a pencil of the form $\pencil(l)$, $l\in\Fn X$, cannot
have $3$-fibers. Hence, all pencils are of type $(0,q)$, $q\le12$, and
$\val l\le12$ for any line $l\in\Fn X$.

\subsection{Locally elliptic configurations}\label{s.locally.elliptic}
A configuration~$S$ or graph $\Fn S$ are called \emph{locally elliptic}
if $\val l\le3$ for
each line $l\in\Fn S$.

\lemma\label{lem.locally.elliptic}
If the graph $\Fn X$ is locally elliptic, then $\ls|\Fn X|\le31$.
\endlemma

\proof
If $\Fn X$ is elliptic, then $\ls|\Fn X|\le\rank\NS(X)\le22$.
If $\Fn X$ contains a plane $a_1,\ldots,a_4$, then $\ls|\Fn X|=4$, as any
other line would increase the valency of one of~$a_i$, see
\autoref{lem.plane}.
Thus, assume that $\Fn X$ is triangle free and contains an affine Dynkin
diagram, which cannot be of type $\tD_4$.
Choose a diagram $\aD\subset\Fn X$ such that
$\mu(\aD)$
is minimal possible and let $\pi:=\pi_{\aD}$, see \autoref{s.trig.free}.

Analyzing affine Dynkin diagrams one-by-one and using the minimality
of~$\aD$, one can easily see that, unless $\aD=\smash{\tD_5}$,
the last term in~\eqref{eq.pencil.count}
is bounded by~$6$.
If the pencil $\pi$ is elliptic, then $\flines\pi\le24$ by
\autoref{cor.count}.
If $\pi$ is quasi-elliptic, \autoref{th.quasi} implies that
$\kchar=2$ and
the parabolic fibers~$F_i$ of~$\pi$ are of type
$\smash{\tD_{2k}}$, $k\ge3$, \smash{$\tE_7$}, or
\smash{$\tE_8$}, with $\sum\mu(F_i)\le20$.
The number of parabolic fibers is at most~$4$ and
$\flines\pi\le24$, see \autoref{cor.count}.
In any case, $\ls|\Fn X|\le30$.

If $\aD=\smash{\tD_5}$, two vertices can be attached to
each of the four monovalent vertices of~$\aD$. However, if a line $s\in\Fn X$
is a section of~$\pi$, then $\pi$ has at most three parabolic fibers,
as otherwise $\val s\ge4$.
As the types of the parabolic fibers are \smash{$\tD_5$} or \smash{$\tA_7$}, with
at least one \smash{$\tD_5$},
it follows that at least one fiber is elliptic; hence,
we have
$\flines\pi\le23$, see~\eqref{eq.Euler}, and $\ls|\Fn X|\le31$.
\endproof

\subsection{Quadrangle-free configurations}\label{s.quad.free}
A configuration $S$ (or graph $\Fn S$)
is said to be \emph{quadrangle free} if $\Fn S$
has no cycles of length~$3$ or~$4$.

\lemma\label{lem.p-q}
Assume that lines $a,b,a_1,\ldots,a_p,b_1,\ldots,b_q\in\Fn X$
\figure
\cpic{p=q}
\caption{The configuration in \autoref{lem.p-q}}\label{fig.p-q}
\endfigure
intersect as shown in \autoref{fig.p-q}, \ie,
\[*
a\cdot b=a\cdot a_i=b\cdot b_i=1,\
a\cdot b_i=a_i\cdot b=a_i\cdot b_j=0\
\]
for all $i=1,\ldots,p$ and $j=1,\ldots,q$.
Then, up to reordering the pair $(p,q)$, one has
either $p,q\le6$, or $p=4$ and $q\le8$, or $p\le3$ and $q\le11$.

Furthermore, in the case $p=q=6$,
if the configuration $\Fn X$ is quadrangle free and
$\kchar\ne2$, one also has $\val a_i\le5$, $\val b_i\le5$ for all
$i=1,\ldots,6$.
\endlemma

\proof
If $p=q=6$, we have a relation
\[
2a+a_1+\ldots+a_6=2b+b_1+\ldots+b_6.
\label{eq.p=q=6}
\]
Hence, $p\ge5$ implies $q\le6$, \cf. \autoref{rem.rels}.
If $p=4$ and $q\ge9$, then
\[*
e:=h+2a+2a_1+2a_2+a_3+a_4-3b-b_1-\ldots-b_9
\]
is an exceptional divisor, see
\autoref{th.Saint-Donat}\iref{Saint-Donat.exceptional}.
Finally, one always has $q\le11$, as the pencil $\pencil(b)$ cannot have more than
twelve fibers, see \autoref{prop.Euler}.

For the last statement, due to relation~\eqref{eq.p=q=6}, it suffices to
show that $\val a_1<7$. By the same relation, any line
$c\in\pencil(a_1)\sminus a$
intersects exactly one of~$b_i$ and is disjoint from all other lines. Assume
that there are six such lines $c_1,\ldots,c_6$, so that $c_i\cdot
b_j=\delta_{ij}$, and consider the sublattice $S\subset\NS(X)$ generated
by~$h$ and all lines $a$, $b$, $a_i$, $b_i$, $c_i$.
We have $\rank S=19$ and $\discr_2S$ is the $\F_2$-vector space generated by
$\frac12(h+a+b+a_1+\ldots+a_6)$ and $\frac12(a_i+a_{i+1})$, $i=2,\ldots,5$.
The isotropic vectors are those of the form
$\frac12(a_2+\ldots+a_6-a_i)$, $i=2,\ldots,6$, and they are all
represented by exceptional divisors, see
\autoref{th.Saint-Donat}\iref{Saint-Donat.exceptional}.
Hence, $S$ has no admissible extensions of index~$2$ and,
by \autoref{th.realization}, $S$ does not admit a geometric realization.
\endproof

\lemma\label{lem.quadrangle.free}
If $\kchar\ne2$ and $\Fn X$ is quadrangle free, then $\ls|\Fn X|\le44$.
\endlemma

\proof
Let $l_0\in\Fn X$ be a line of maximal valency. In view of
\autoref{lem.locally.elliptic}, we can assume that $\val l_0\ge4$. Choose four
lines $l_1,\ldots,l_4\in\Fn X$ adjacent to~$l_0$; together with~$l_0$, they
constitute a type~$\tD_4$ affine Dynkin diagram contained in $\Fn X$.
Applying~\eqref{eq.pencil.count} to the pencil $\pi:=\pi_{\aD}\:X\to\Cp1$,
we have
\[*
\ls|\Fn X|=\flines\pi+\val l_0+\val l_1+\ldots+\val l_4-8.
\]

The valencies are estimated by applying
\autoref{lem.p-q} to $a=l_i$, $i=1,\ldots,4$, and $b=l_0$
(recall that we assume $\val l_0\ge\val l_i$).
In the worst case, where $\val l_i=6$, $i=0,\ldots,4$, we have $\sum_i\val l_i=30$;
otherwise, $\sum_i\val l_i\le29$.

The parabolic fibers of~$\pi$ are of types $\tD_4$ or $\tA_5$,
with at least one $\tD_4$, and its elliptic fibers are of linear types
$\bA_p$, $p\le4$, or $\bA_1^2$. Thus, we have $\flines\pi\le23$, see
\autoref{cor.count}. We assert that, if $\val l_0=6$, the pencil has at
most three parabolic fibers and, hence, $\flines\pi\le22$;
these inequalities imply $\ls|\Fn X|\le44$ in the statement.
For the last assertion, let $n_1,n_2\ne l_1,\ldots,l_4$ be
the two extra lines adjacent
to~$l_0$. They are bisections of~$\pi$. If $F$ is another type~\smash{$\tD_4$}
fiber, then either one of $n_1$, $n_2$ intersects two monovalent vertices of~$F$
or both lines intersect the central vertex; in either case, we obtain a
quadrangle. If $F_1$, $F_2$, $F_3$ are three type~\smash{$\tA_5$} fibers,
then, each~$n_i$ intersecting two lines in each~$F_j$, we obtain
$\val n_i\ge7$, which contradicts to our choice of~$l_0$.
\endproof

\subsection{Configurations with a quadrangle}\label{s.quad}
Assume that a triangle free configuration $\Fn X$ contains a quadrangle
$Q:=\{l_1,l_2,l_3,l_4\}$. (The lines constituting a quadrangle
are always listed according to their cyclic order in the affine Dynkin
diagram).
Let $\pi_Q\:X\to\Cp1$ be the corresponding elliptic pencil.
It has a certain number~$s_3$ of parabolic fibers of type~\smash{$\tA_3$} and,
for $p=1,2$, a certain number~$s_p$ of elliptic fibers of
linear type~$\bA_p$.
By \autoref{cor.count}, we have
\[
\flines{\pi_Q}=4s_3+2s_2+s_1,\qquad
4s_3+3s_2+2s_1\le24.
\label{eq.quad.counts}
\]
Identity~\eqref{eq.pencil.count} becomes
\[
\ls|\Fn X|=\flines{\pi_Q}
 +\val l_1+\ldots+\val l_4-\dif(l_1,l_3)-\dif(l_2,l_4)-8,
\label{eq.quad}
\]
where the correction terms
$\dif(l_i,l_j):=\ls|\pencil(l_i)\cap\pencil(l_j)|-2$
are nonnegative.

\lemma[see~\cite{DIS}]\label{lem.quad<=20}
If $\kchar\ne2$ and $\Fn X$ is triangle free, then $\flines{\pi_Q}\le21$.
\endlemma

\proof
In view of~\eqref{eq.quad.counts},
we only need to eliminate the triples $(s_3,s_2,s_1)=(6,0,0)$, $(5,1,0)$,
and $(5,0,2)$, the former being a consequence of the two latter.
To this end, we consider the lattice~$S$ generated by~$h$ and the lines in
the fibers,
use an analog of \autoref{lem.discr} to estimate from below the $2$-torsion
of the discriminants of admissible finite index extensions $\tS\supset S$
such that $\Fn\tS$ is still triangle free, and
apply \autoref{th.realization} (the part related to $p=2$)
to show that $\tS$ does not admit a primitive
geometric realization.
Details are left to the reader.
\endproof

Now, the following lemma is an immediate consequence of~\eqref{eq.quad},
\autoref{lem.quad<=20}, and
the bound $\val l\le12$ for each line $l\in\Fn X$.

\lemma\label{lem.quad}
If $\kchar\ne2$ and $\Fn X$ is triangle free and
contains a quadrangle, then $\ls|\Fn X|\le61$.
\endlemma

\subsection{Proof of \autoref{th.triangle.free}}\label{proof.triangle.free}
If $X$ is not supersingular, we can use \autoref{th.char=0} and assume that
$\kchar=0$. Hence, depending on the type of the configuration $\Fn X$, the
statement of the theorem follows from
\autoref{lem.quad} ($\Fn X$ contains a quadrangle),
\autoref{lem.quadrangle.free} ($\Fn X$ is quadrangle free and has a line of
valency at least~$4$) and \autoref{lem.locally.elliptic}
($\val l\le3$ for all $l\in\Fn X$).
\qed

\section{Exotic pencils}\label{S.exotic}

\emph{Exotic} are pencils contained in a supersingular quartic over an
algebraically closed field
of characteristic~$3$. The most interesting feature of such pencils is the
fact that the existence of a pencil of a certain type does not guarantee the
existence of a maximal pencil of any smaller type, see
\autoref{cor.exotic.type} below.

\subsection{Quasi-elliptic pencils}\label{s.quasi}
According to \autoref{prop.Euler}, a pencil
of the form $\pi[l]$, $l\in\Fn X$, is quasi-elliptic if and only if
$\kchar=3$ and $\pencil(l)$ is of type $(10,0)$.
Since $\rank P_{10,0}=22$, the lattice $\NS(X)$ is a
finite index extension of~$P_{10,0}$; its kernel is denoted by
$\CX:=\NS(X)/P_{10,0}\subset\discr P_{10,0}$. Note that $\CX$
is an $\F_3$-vector space, as so is $\discr P_{10,0}$, see
\autoref{obs.discr}.

\proposition\label{prop.(10,0)}
The map $s\mapsto(s\bmod P_{10,0})\in\CX$ establishes a bijection between the
set of sections of $\pencil(l)$ and the set
$\{\Ga\in\CX\,|\,\Ga\cdot\Gl=\frac23\bmod\Z\}=\CX\cap\oO_0$.
\endproposition

\proof
Since $\NS(X)\subset P_{10,0}\dual$, each section represents an element
$\Ga\in\CX$ as in the statement.
Distinct sections represent distinct elements by \autoref{lem.10-2} and,
by the same lemma, any element in the
$\OG_h(P_{10,0})$-orbit of $\Go+6\Gl$ is a section.
\endproof

\remark\label{rem.(10,0)}
If $s_1\ne s_2$ are two sections, then
$s_1-s_2\bmod P_{10,0}$ is an isotropic vector orthogonal to~$\Gl$.
By \autoref{obs.discr} and \autoref{lem.discr}, this difference is in
$\oM_6\cup\oM_9$.
Each section intersects all ten fibers of~$\pencil(l)$, see
\autoref{lem.plane}; hence, applying \autoref{lem.10-2}
(see also \autoref{rem.rels})
to~$l$, $s_1$, and the ten common lines, we conclude that
$s_1\cdot s_2=1$ if and only if $s_1-s_2\in\oM_9$.
\endremark

\corollary\label{cor.(10,0)}
Let $\pencil\subset\Fn X$ be a pencil of type $(10,0)$. Then
\roster*
\item
$\sect\pencil=0$ or $3^r$, $0\le r\le4$\rom;
respectively, $\ls|\Fn X|=31$ or $31+3^r$ and, in the latter case,
one has $\Artin(X)=5-r$.
\endroster
There is a unique configuration of size~$112$ and
two configurations of size~$58$\rom; they are as in
\autoref{th.char=3}\iref{112.(10,0)} and \iref{58.(10,0)x2},
\iref{58.(10,0)}, respectively.
\endcorollary

\proof
According to \autoref{prop.(10,0)}, the sections of~$\pencil$
constitute the affine subspace $\CX\cap\oO_0$; it is either empty
or has dimension $r:=\dim\CX-1\le4$. In the latter case, up to automorphism,
we have
$\CX=\F_3\Go\oplus\CX'$, where
$\CX':=\CX\cap\Go^\perp\cong(\CX\cap\Gl^\perp)/\Gl$ is a ternary code of
length $10$ with all Hamming norms~$6$ or~$9$ (see \autoref{lem.discr}).
The number of such codes of dimension~$1$, $2$, $3$, $4$ is,
respectively, $2$, $3$, $2$, $1$.
(The statement on codes of dimension~$4$, \ie, the uniqueness of the
Fermat quartic, can also be established by other means, see
\autoref{rem.Fermat} below.)
\endproof

We refer to \autoref{s.eq.(10,0)} for a geometric description of
quasi-elliptic pencils.

\subsection{The lattice $\NS(X)$}\label{s.NS.char=3}
We often need to establish the
\hbox{(non-)}\allowbreak uniqueness in the genus
of the orthogonal complement of $P_{p,q}$ in $\NS(X)$.
To do so, we use the results of \autoref{s.discr.lemmas} and start with a
``minimal'' lattice~$T$, \ie, the one with largest discriminant.
More precisely, \autoref{prop.genus.T} and \autoref{obs.discr} imply that
\[*
\discr_2T\cong\<\tfrac12\>^q,\qquad
\discr_rT=0\ \ \text{for $r>3$}.
\]
The $9$-torsion is also determined
by~$(p,q)$: in the terminology of \autoref{s.discr.lemmas}, we have
\[*
(3\discr_3T,\qc)\cong(3\discr_3P_{p,q},-\qc).
\]
However, the $3$-elementary part of $\discr_3T$ is not fixed,
and we make it as large as possible.
Then,
$\NS(X)$ is a finite index extension of $N:=P_{p,q}\oplus T$, and it this
extension that we try to describe by means of its kernel
\[*
\CX:=\NS(X)/N\subset\discr N=\discr P_{p,q}\oplus\discr T.
\]
We reserve this notation till the end of the section.
Clearly, $\CX=\bigoplus_r\CX_r$, where we let $\CX_r:=\CX\cap\discr_rN$ for a
prime~$r$.
Then, $\CX_r=0$ for $r>3$ and $\CX_2$ is the graph of a bijective anti-isometry
\[*
\psi_2\:\discr_2P_{p,q}\longto\discr_2T.
\]
(The case where the quotient $\NS(X)/P_{p,q}$ has $2$-torsion is discussed
separately in \autoref{s.(0,q)}.)
Our principal concern is the natural map
\[*
\sec_3\:\{\text{sections of $\pencil$}\}\longto\CX_3
\]
sending a section~$s$ to the projection of $(s\bmod N)$ to $\CX_3$. All
nonempty fibers of the composition of~$\sec_3$ and the projection to
$\discr_3P_{p,q}$ are over the affine space
\[*
\oO_*:=\bigl\{\Ga\in\discr_3P_{p,q}\bigm|\Ga\cdot\Gl=-\tfrac13r\bigr\}=
 \oO_0\cup\oO_+\cup\oO_-.
\]
Since we do not always assume a pencil $\pencil\subset\pencil(l)$
maximal, there is a similar map
\[*
\fib_3\:\pencil(l)\sminus\pencil\longto\CX_3
\]
with all nonempty fibers over~$-\oL\cup\oZ$.

We make a few general observations. Fix an anti-isometry~$\psi_2$ as above;
it defines an extension $N_3\supset N$ with $\discr_rN_3=0$ for $r\ne3$.
Moreover, if $p+q=1\bmod3$ (which we usually assume), the lattice $N_3$ is
$3$-elementary. Recall that the image of $\OG_h(P_{p,q})$ in
$\Aut\discr_2P_{p,q}$ is the symmetric group $\SG{q}$, see
\autoref{obs.discr}.
Let
\[*
\OG(T,\psi_2):=\bigl\{g\in\OG(T)\bigm|
 \psi_2\1\circ\bar g\circ\psi_2\in\SG{q}\bigr\},
\]
where $\bar g$ is the image of~$g$ in $\Aut\discr_2T$. Then, since the
actions of $\OG_h(P_{p,q})$ on $\discr_2P_{p,q}$ and $\discr_3P_{p,q}$ are
independent, see \autoref{obs.discr}, the action of
the group $\OG_h(N_3)$ on
$\discr N_3$ reduces to the product action of
$\OG_h(P_{p,q})\times\OG(T,\psi_2)$, so that the orbits of elements of
order~$3$ are products of the form
\[*
\oM_k\times\ooO,\ k=0,\ldots,p,\quad
\pm\oL\times\ooO,\quad
\pm\oO_s\times\ooO,\ s\in\F_3,
\]
where $\ooO\subset\discr_3T$ is an $\OG(T,\psi_2)$-orbit.
When describing the fibers of the maps $\sec_3$ and $\fib_3$ and discussing
the admissibility of the extension (the non-existence of exceptional
divisors), it suffices to check a single representative of each orbit.

Another observation concerns the choice of
the anti-isometry~$\psi_2$. We have the following
obvious lemma.

\lemma\label{lem.psi2}
Assume that a generator $3\Gn_k\in\discr_2P_{p,q}$, $k=1,\ldots,q$, is such
that the image $\psi_2(3\Gn_k)$ is represented by $\frac12a$, where $a\in T$,
$a^2=-6$. Then, the vectors
\[*
n_k+3\Gn_k\pm\tfrac12a=-\Gl+\Gn_k\pm\tfrac12a\in\NS(X)
\]
are lines that belong to the maximal pencil $\pencil(l)$
containing~$\pencil$. As a consequence,
the $k$-th $1$-fiber of~$\pencil$ becomes a $3$-fiber of $\pencil(l)$ and
$\pencil$ is not maximal.
\endlemma

\subsection{Pencils of type $(7,0)$}\label{s.(7,0)}
By \autoref{prop.types}, any pencil of type $(p,q)=(7,0)$ is exotic, and
the minimal lattice~$T$ is given by
\autoref{lem.3-discr} as $\bbE_6:=\bE_6\dual(3)$.
The homomorphism
$\OG(\bbE_6,\psi_2)=\OG(\bbE_6)=\OG(\bE_6)\to\Aut\discr\bbE_6$ is an
isomorphism;
hence, the $\OG(\bbE_6)$-orbits
in $\discr\bbE_6$ are $\oE_0:=\{0\}$ and
\[*
\textstyle
\oE_r:=\bigl\{\Ga\in\discr\bbE_6\bigm|
 \Ga\ne0,\ \Ga^2=-\frac23r\bmod2\Z\bigr\},\quad r=1,2,3,
\]
and each $\Ga\in\oE_r$
has
a representative of the form $\frac13a$,
where $a\in\bbE_6$, $a^2=-6r$.
The following statement is immediate, \cf. \autoref{lem.discr} or
\autoref{prop.(10,0)}: one can easily check all orbits one by one.

\proposition\label{prop.(7,0)}
Let $\pencil\subset\Fn X$ be a pencil of type $(7,0)$. Then
the kernel~$\CX$ is disjoint from the $\OG_h(N)$-orbits
\[*
\oZ\times\oE_3,\quad
\oM_1\times\oE_2,\quad \oM_2\times\oE_1,\quad \oM_3\times\oE_0,\quad
\pm\oL\times\oE_0,\quad \pm\oO_0\times\oE_0.
\]
The map $\sec_3$ is a bijection onto $\CX\cap(\oO_1\times\oE_1)$.
The pencil~$\pencil$ is not
maximal if and only if it extends to a quasi-elliptic pencil if and only if
$\CX\cap(\oL\times\oE_3)\ne\varnothing$.
\endproposition

\corollary\label{cor.(7,0)}
Let $\pencil\subset\Fn X$ be a maximal pencil of type $(7,0)$.
Then\rom:
\roster*
\item
$\sect\pencil=36$ or $\sect\pencil\le27$\rom; respectively,
$\ls|\Fn X|=58$ or $\ls|\Fn X|\le49$\rom;
\item
if $\sect\pencil>5$, then $\sect\pencil=0\bmod3$.
\endroster
These bounds are sharp, and there is a unique configuration
$\Fn X$ of size~$58$\rom; it is as in \autoref{th.char=3}\iref{58.(7,0)}.
\endcorollary

\proof
By \autoref{prop.(7,0)}, we have $\CX\cap\discr\bbE_6=0$;
hence, $\CX$ is the graph of a certain anti-isometry
$\psi\:\DD\to\discr\bbE_6$, where
the domain $\DD\subset\discr P_{7,0}$ is disjoint
from~$\oL$, $\oM_1$, $\oM_2$.
The projection $\DD\cap\Gl^\perp\to\Gl^\perp\!/\Gl$ is injective
and its image
is a ternary code of length~$7$ and
minimal Hamming distance~$3$; it has dimension at most~$4$.
Hence, $\dim\DD\le5$.
The sections of~$\pencil$ are in a one-to-one correspondence with the
vectors $\Ga\in\DD':=\DD\cap\oO_*$
such that $\psi(\Ga)\in\oE_1$, \ie, $\Ga^2=\frac23$.
Thus, the congruence in the statement follows from \autoref{lem.counts}.

If $\dim\DD\le4$, then $\sect\pencil\le\ls|\DD'|\le27$;
henceforth, we assume that $\dim\DD=5$.

If $\Ker\psi=0$, then $\psi$ restricts to an injective map from~$\DD'$ to a
proper affine subspace $\CE'\subset\discr\bbE_6$ disjoint from~$0$.
For any such space, $\ls|\CE'\cap\oE_1|\le27$.

Assume that $\Ker\psi\ne0$. By \autoref{prop.(7,0)} (or \autoref{lem.discr}),
$\dim\Ker\psi=1$ and $\Ker\psi$ is generated by an element of~$\oM_6$. In
this case, $\psi$ restricts to a three-to-one map $\DD'\onto\CE'$, where
$\CE'\subset\discr\bbE_6$ is an affine subspace disjoint
from~$0$ and $\dim\CE'=3$. With one exception, one has $\ls|\CE'\cap\oE_1|\le9$.
In the exceptional case, $\ls|\CE'\cap\oE_1|=12$ (hence $\sect\pencil=36$),
both spaces $\CE_0:=\psi(\Gl^\perp)\subset\CE:=\Im\psi$ are
nondegenerate, $\Br\CE_0=2$, and $\Br\CE=0$; in other words,
\[*
\CE_0\cong\bigl\<\tfrac43\bigr\>^3\subset
 \CE=\CE_0\oplus\bigl\<\tfrac43\bigr\>\subset
 \discr\bbE_6=\CE\oplus\bigl\<\tfrac23\bigr\>.
\]
On the other hand, in the space
$\discr_3P_{7,0}$, there is a unique pair $(\DD,\DD\cap\Gl^\perp)$ satisfying
all the assumptions above and anti-isomorphic (modulo kernel) to $(\CE,\CE_0)$.
Since the stabilizer of $(\CE,\CE_0)$ restricts to the full
group $\Aut\CE_0$,
the anti-isometry $\psi\:\DD\to\CE$ is also unique; it gives rise to a
quartic as in the statement.
\endproof

\subsection{Exotic pencils of type $(4,6)$}\label{s.(4,6)}
Consider a pencil $\pencil$ of type~$(4,6)$ and \emph{assume} that it is
exotic. By \autoref{lem.2-3-discr},
the minimal orthogonal complement~$T$ of~$P_{4,6}$
is $\bH_6(6)$,
which is unique in its genus. The $\OG(T)$-orbits in $\discr_3T$ are
\[*
\oH_r:=\bigl\{\Ga\in\discr_3T\bigm|\|\Ga\|=r\bigr\},\quad r=0,\ldots,6,
\]
where $\|\cdot\|$ is the Hamming norm in the obvious basis.

By \autoref{lem.Dima}, there are two
essentially distinct choices for the anti-isometry
$\psi_2\:\discr_2P_{4,6}\to\discr_2T$, \viz. the identity and $\rt_1$.
In the former case, $\pencil$ extends to a quasi-elliptic pencil,
see \autoref{lem.psi2}.
Thus, from now on we assume that $\psi_2=\rt_1$;
such a pencil $\pencil\subset\NS(X)$ is called \emph{non-trivially} exotic.
We have $\OG(T,\psi_2)=\OG(T)$, and, as above, the following statement is
straightforward.

\proposition\label{prop.(4,6)}
Let $\pencil\subset\Fn X$ be a non-trivially exotic pencil of
type $(4,6)$. Then the kernel~$\CX$ is disjoint from the $\OG_h(N)$-orbits
\[*
\oZ\times\oH_3,\quad
\oM_1\times\oH_*,\quad
\oM_2\times\oH_1,\quad
\oM_3\times\oH_0,\quad
\pm\oL\times\oH_*,\quad
\pm\oO_1\times\oH_1
\]
\rom(where $*$ stands for whichever index is appropriate\rom).
Furthermore, the map $\sec_3$ is
\roster*
\item
three-to-one over $\CX\cap(\oO_2\times\oH_2)$ and
\item
one-to-one over $\CX\cap(\oO_0\times\oH_*)$ and
$\CX\cap(\oO_1\times\oH_4)$\rom;
\endroster
all other fibers are empty.
The pencil~$\pencil$ is not
maximal if and only if it extends to a quasi-elliptic pencil
if and only if
$\CX\cap(\oZ\times\oH_6)\ne\varnothing$.
\endproposition

\corollary\label{cor.(4,6)}
Let $\pencil\subset\Fn X$ be a maximal exotic pencil of type $(4,6)$.
Then\rom:
\roster*
\item
$\sect\pencil=39$ or $\sect\pencil\le33$\rom; respectively,
$\ls|\Fn X|=58$ or $\ls|\Fn X|\le52$\rom;
\item
if $\sect\pencil>19$, then $\sect\pencil=0\bmod3$.
\endroster
These bounds are sharp. If $\ls|\Fn X|=58$, then
$\Fn X$ contains a pencil of type~$(10,0)$ or $(7,0)$ and, hence,
is as in \autoref{th.char=3}\iref{58.(10,0)} or~\iref{58.(7,0)}.
\endcorollary

\proof
Consider the projections $\DD\subset\discr_3P_{4,6}$ and
$\CH\subset\discr_3T$ of the kernel~$\CX_3$ to the two summands of $\discr N$. Since
$\pencil$ is assumed maximal, \autoref{prop.(4,6)} implies that
$\CX_3\cap\discr_3T=0$;
hence, $\CX_3$ is the graph of a certain anti-isometry $\psi\:\DD\onto\CH$.
Conversely,
an anti-isometry
$\psi\:\DD\onto\CH$ gives rise
to a nonsingular quartic if and only if
\roster
\item\label{4-6.U}
$\DD$ is disjoint from the orbits $\oM_1$ and $\pm\oL$,
\item\label{4-6.T}
$\CH$ is disjoint from the orbit $\oH_1$; in other words,
$\CH\subset\discr_3T$
is a ternary code of minimal Hamming distance at least~$2$, and
\item\label{4-6.ker}
$\Ker\psi$ is disjoint from the orbit~$\oM_3$.
\endroster
The first condition implies that $\dim\DD\le4$ and, up to the action of
$\OG_h(P_{4,6})$, there are three subspaces of dimension~$4$, \viz.
\[*
\DD_r:=\bigl[\F_3(\Gm_1+\ldots+\Gm_4)\oplus\F_3(\Go+r\Gl)\bigr]^\perp,\quad
r\in\F_3.
\]

Let $\DD':=\DD\cap\oO_*$.
Since $\Ker\psi\cap\Gl^\perp=0$, the restriction $\psi|_{\DD'}$
is a one-to-one map onto an
affine subspace $\CH'\subset\CH$. Denote $h_r:=\ls|\CH'\cap\oH_r|$,
$r=0,\ldots,6$; we have $h_1=0$ and
$\sect\pencil$ is given by \autoref{prop.(4,6)} as
\[
\sect\pencil=h_0+3h_2+h_3+h_4+h_6.
\label{eq.4-6.count}
\]

If $\dim\DD\le3$, then $\dim\DD'\le2$ and, even
if each element of~$\DD'$ triples,
we have $\sect\pencil\le27$.
Therefore, we consider the three spaces $\DD_r$, $r\in\F_3$.

Let $\DD=\DD_r$ and $r\ne0$. Then $\Go_r:=r\Go-\Gl\in\DD$ and $\DD'$ is
given intrinsically as
\[*
\DD'=\bigl\{\Ga\in\DD\bigm|\Ga\cdot\Go_r=-\tfrac13\bigr\}.
\]
Specializing~$r$ to $\pm1$, we have $\Br\DD_r=2+2r$ and
$\Go_r^2=\frac23r$; since $\Ker\psi=0$, these invariants determine the
isomorphism type of $(\CH,\Gb):=\psi(\DD_r,\Go_r)$. Then, for each
$\OG(T)$-orbit of pairs $(\CH,\Gb)$, the number of sections $\sect\pencil$ is
given by~\eqref{eq.4-6.count}. With condition~\iref{4-6.T} taken into
account, we have
\roster*
\item
six orbits, if $r=1$, and then $\sect\pencil=39,33,27,27,24,21$, or
\item
three orbits, if $r=-1$, and then
$\sect\pencil=24,21,18$.
\endroster

Let $\DD=\DD_0\cong(\DD_0/\ker)\oplus\F_3\Go$, $\Go^2=0$.
In this case, we have $\Br(\DD/\ker)=2$ and
$\DD'=(\DD/\ker)+\Go$. If $\Ker\psi=0$, the pair $(\CH,\CH')$ is
anti-isomorphic to $(\DD,\DD')$;
otherwise, $\psi$ maps
$\DD'$ onto $\CH'=\CH$, which
is anti-isomorphic to $\DD/\ker$. Taking condition~\iref{4-6.T} into
account and using~\eqref{eq.4-6.count}, we obtain
\roster*
\item
four orbits, if $\Ker\psi=0$, and then $\sect\pencil=30,27,24,21$, or
\item
four orbits, if $\Ker\psi\ne0$, and then $\sect\pencil=39,33,27,21$.
\endroster

There are two cases where $\sect\pencil=39$.
If $r\ne0$, the stabilizer of $(\DD_r,\Go_r)$ induces the full
group $\Aut(\DD_r,\Go_r)$.
Similarly, the stabilizer of $\DD_0$ induces the full
group $\Aut(\DD_0/\ker)$. Hence, in both cases,
an anti-isometry
$\psi\:\DD\onto\CH$ is unique up to the two-sided action of
the group
$\OG_h(P_{4,6})\times\OG(\bH_6)$.
Choosing a representative and computing all lines, we
conclude that the quartic is one of those considered above as it
contains a pencil (other than~$\pencil$)
of type $(10,0)$ or $(7,0)$.

For the congruence, rewrite~\eqref{eq.4-6.count} in the form
\[*
\sect\pencil=\ls|\CH'|+3h_2-(h_2+h_5).
\]
In the notation of \autoref{lem.counts}, we have
$h_2+h_5=n_2(\CH')$. By the lemma,
either $\ls|\CH'|\le9$ and $h_2\le h_2+h_5\le5$, and then
$\sect\pencil\le19$, or $\sect\pencil=0\bmod3$.
\endproof

\subsection{Exotic pencils of type $(4,0)$}\label{s.(4,0)}
Consider a pencil $\pencil$ of type~$(4,0)$ and assume that it is
exotic.
For the lattice $\barT$ in \autoref{lem.3-discr},
we have $\discr\barT=\<\frac23\>^2$; hence, $\barT$ is
the orthogonal complement of a sublattice $\bA_2^2\subset U$, where
$U$ is even unimodular of rank~$16$,
\ie, $U=\bE_8^2$ or $\bD_{16}^+$ (the unique,
up to isomorphism, even unimodular extension of $\bD_{16}$).
Then, the minimal orthogonal
complement
is $T=\barT\dual\!(3)$.

The sublattice $\bA_2^2$ is embedded into the maximal root system contained
in~$U$; there are two embeddings to~$\bE_8^2$ and one embedding to $\bD_{16}$.

If $U=\bE_8^2$, we have $T=\bbE_6^2$ or $\bA_2^2\oplus\bE_8$
(as obviously $\bA_2\dual(3)\cong\bA_2$), and
the latter
lattice contains $(-2)$-vectors, contradicting
\autoref{th.Saint-Donat}\iref{Saint-Donat.exceptional}.

If $U=\bD_{16}^+$ and $\bD_{16}$ is represented as
the maximal even sublattice in $\bH_{16}$,
so that $\bD_{16}^+$ is the extension by $\frac12\be$,
we can assume that the two copies of $\bA_2\subset\bD_{16}$ are generated by
the pairs of roots
$\{\e_1-\e_2,\e_2-\e_3\}$ and $\{\e_4-\e_5,\e_5-\e_6\}$; then,
$\frac13(\e_1+\ldots+\e_6)\in T$ is a $(-2)$-vector.

Thus, we are left
with $T=\bbE_6^2$. The orbits of the $\OG(T)$-action on $\discr T$
are described in \autoref{s.(7,0)}; we will abbreviate
$\oE_{r,s}:=(\oE_r\times\oE_s)\cup(\oE_s\times\oE_r)$.

\proposition\label{prop.(4,0)}
Let $\pencil\subset\Fn X$ be an exotic pencil of type $(4,0)$. Then
\roster*
\item
$\CX\cap\discr P_{4,0}=0$ and
$\CX\cap\discr T\subset\oE_{0,0}\cup\oE_{3,3}$\rom;
\item
$\CX$ is disjoint from
$\oM_1\times\oE_{1,1}$, $\oM_1\times\oE_{0,2}$, $\oM_2\times\oE_{0,1}$,
and $\oO_1\times\oE_{0,1}$.
\endroster
Furthermore, the map $\sec_3$ is
\roster*
\item
three-to-one over $\CX\cap(\oO_2\times\oE_{0,2})$ and
\item
one-to-one over $\CX\cap(\oO_2\times\oE_{1,1})$\rom;
\endroster
all other fibers are empty.
The pencil~$\pencil$ is not maximal if and only if $\CX$ intersects one of
the following orbits\rom:
\roster*
\item
$\oL\times\oE_{0,3}$,
and then $\pencil$ is contained in a pencil of type $(7,0)$, or
\item
$\oL\times\oE_{1,2}$,
and then $\pencil$ is contained in a pencil of type $(4,3)$.
\endroster
\endproposition

As an immediate consequence, for a maximal exotic pencil of type $(p,q)$,
$p\ge4$, we
have either $p=4$ and $q=0$, $3$, $6$ or $(p,q)=(7,0)$ or $(10,0)$.

\subsection{Exotic pencils of type $(0,q)$}\label{s.(0,q)}
Consider an
exotic pencil~$\pencil$ of type $(0,10)$. Here, the situation is much more
diverse than in the previous sections.

First, by \autoref{lem.2-3-discr}, the
minimal orthogonal complement~$T$ has the form $\barT(6)$, where $\barT$ is
an odd unimodular lattice of rank~$10$, \ie, $\barT=\bH_{10}$ or
$\bE_8\oplus\bH_2$.

Second, the intersection $\CX_2\cap\discr_2P_{0,10}$ may be nontrivial:
it may contain an element
$\Ga\in\oN_8\subset\discr_2P_{0,10}$ (see \autoref{lem.discr}), and
then also
$\CX_2\cap\discr_2T\ne0$.
If this is the case, we still
choose an anti-isometry $\psi_2\:\discr_2P_{0,10}\to\discr_2T$
and represent~$\CX_2$ as
$\F_2\Ga\oplus\operatorname{graph}(\psi_2|_{\Ga^\perp}\!)$.
If $\psi_2$ is fixed, the orbit~$\oN_8$ splits into the orbits of the
subgroup
\[*
G:=\bigl\{g\in\SG{q}\bigm|
 \text{$\psi_2\circ g\circ\psi_2\1$ is in the image of $\OG(T)$}\bigr\},
\]
and it suffices to consider for~$\Ga$ a single representative of each suborbit.
Then, with $\Ga\ne0$ fixed, the group $\OG(T,\psi_2)$ in
\autoref{s.NS.char=3} can be extended to the larger subgroup
\[*
\OG(T,\psi_2,\Ga):=\bigl\{g\in\OG(T)\bigm|
 \text{$(\psi_2\1\circ\bar g\circ\psi_2)|_{\Ga^\perp}\!=s|_{\Ga^\perp}$
 for some $s\in\SG{q}(\Ga)$}\bigr\},
\]
where $\SG{q}(\Ga)$ is the stabilizer of~$\Ga$.
\autoref{lem.psi2} should be restricted to the generators
$\Gn_k\in\Ga^\perp$, leaving more choice for~$\psi_2$.
Other observations made in \autoref{s.NS.char=3} apply literally.

Finally, there are several choices for~$\psi_2$ itself:
if $\barT=\bH_{10}$, the nine classes are given by \autoref{lem.Dima}, and
if $\barT=\bE_8\oplus\bH_2$, the inverse~$\psi_2\1$ is determined by the
Hamming norms of the images of the two generators
$\frac12\e_1,\frac12\e_2\in\discr_2\bH_2(6)$,
which may be $(1,9)$ or $(5,5)$. (Recall that the
characteristic vector
$\frac12(\e_1+\e_2)$ is mapped to the characteristic vector
and the map $\OG(\bE_8)\to\Aut\discr_2\bE_8(6)$ is surjective.)

Summarizing, we have two choices for~$\barT$, $(9+2)$ choices for~$\psi_2$,
and, for each~$\psi_2$, up to four choices for $\Ga\in\oN_0\cup\oN_8$.
Furthermore, in most cases, there is no simple description of the
group $\OG(T,\psi_2,\Ga)$ and its
orbits on $\discr_3T$.
Hence, we use \GAP~\cite{GAP4} to enumerate the orbits and, afterwards,
quartics. (We disregard the quartics in which a $1$-fiber
of~$\pencil$ becomes a $3$-fiber;
however, we do allow extra $1$- and $3$-fibers.)
The resulting statement is as follows.

\proposition\label{prop.(0,q)}
Let $\pencil\subset\Fn X$ be a maximal exotic pencil of type $(0,q)$,
$q\ge10$.
\roster*
\item
If $q=10$, then $\sect\pencil\le29$\rom; hence, $\ls|\Fn X|\le40$.
\item
If $q=11$, then $\sect\pencil\le22$\rom; hence, $\ls|\Fn X|\le34$.
\item
If $q=12$, then $\sect\pencil\le39$\rom; hence, $\ls|\Fn X|\le52$.
\endroster
If $\Fn X$ as above is triangle free, then $\ls|\Fn X|\le37$\rom;
this bound is sharp and it is attained by a unique configuration\rom:
\roster*
\item
$\ps=\PS{[ [ 0, 12 ], 1 ], [ [ 0, 9 ], 22 ], [ [ 0, 6 ], 14 ]}$,
$\Artin(X)=2$, and $\rank\Fano(X)=21$.
\endroster
\endproposition

\remark\label{rem.char=3.counts}
The counts $\ls|\Fn X|$ observed in the proof of \autoref{prop.(0,q)} are
\[*
\bigl\{11, 12, \ldots, 33, 34, 36, 37, 40, 43, 46, 52 \bigr\}.
\]
Together with Corollaries~\ref{cor.(10,0)}, \ref{cor.(7,0)},
and~\ref{cor.(4,6)}, this list substantiates \autoref{conj.char=3}.
Note that the value $\ls|\Fn X|=49$ is also taken,
as follows from \autoref{prop.(7,0)}
or the proof of \autoref{prop.(4,6)}.
\endremark


\subsection{Types of exotic pencils}\label{s.types.exotic}
In conclusion, we list the possible types of exotic pencils. Comparison of
Tables~\ref{tab.Euler} and~\ref{tab.exotic} explains the term ``exotic.''

\proposition\label{cor.exotic.type}
The type $(p,q)$ of a maximal exotic pencil $\pencil\subset\Fn X$ takes
only the values listed in \autoref{tab.exotic}.
\table
\caption{The types $(p,q)$ of exotic pencils, see \autoref{cor.exotic.type}}\label{tab.exotic}
\centerline{\setbox0\hbox{$0,3,6$}\vbox{\halign{
 \hss\strut$#$&&\quad\hss$#$\hss\cr
p:&10&7&4&3&2&1&0\cr
q:&0&0&0,3,6&\le6&\le8&\le9&\le12\cr
}}}
\endtable
\endproposition

There are a few further restrictions,
which we do not discuss. For example,
a maximal exotic  pencil
of type $(3,q)$ exists if and only if $q\le4$ or $q=6$.
The upper bounds for~$q$ in \autoref{tab.exotic} are sharp.


\proof[Proof of \autoref{cor.exotic.type}]
In view of Propositions~\ref{prop.Euler} and~\ref{prop.(4,0)}, there only
remains
to show that maximal pencils of types $(p,q)=(3,7)$, $(2,9)$, or $(1,10)$ do not
exist.

For the first two types, we start with a subpencil~$\pencil$ of type $(2,8)$.
According to \autoref{lem.discr} (\cf. also \autoref{prop.types.2}), this subpencil
can be chosen so that the quotient group $\NS(X)/P_{2,8}$ does not have
$2$-torsion. Then, the minimal orthogonal complement of~$P_{2,8}$ is
$T=\bH_8(6)$ (see \autoref{lem.2-3-discr}), and, by \autoref{lem.Dima}, there
are three essentially distinct choices for the anti-isometry
$\psi_2\:\discr_2P_{2,8}\to\discr_2T$:
the identity, $\rt_1$, and $\rt_1\rt_3$.
The first two give rise to ``immediate'' extra fibers and embed~$\pencil$ to
a pencil of type
$(10,0)$ or
$(4,6)$, see \autoref{lem.psi2}; hence, we choose the last one.
The resulting quartic is nonsingular; in particular,
we conclude that maximal exotic pencils of type $(2,8)$ do exist.

The orbits of the $\OG(T,\psi_2)$-action on
$\discr_3T$ are characterized by the ``triple'' Hamming norm, \ie, the
sequence $(u,v,w)$ of
Hamming norms in the coordinates $\{1,2\}$, $\{3,4,5,6\}$, and $\{7,8\}$;
we denote these orbits by $\oH_{u,v,w}$. (In fact,
$\OG(T,\psi_2)$ contains also the permutation
$(1,3)(2,4)(5,7)(6,8)$ of the basis vectors,
which we ignore to simplify the description of the orbits.)
As in the previous sections,
we check that the kernel $\CX_3\subset\discr_3N$ defining
the quartic is disjoint from $\pm\oL\times\oH_*$ and
the pencil~$\pencil$ is not maximal if and only if $\CX_3$ intersects any of
\[*
\oZ\times\oH_{2,2,2},\quad
\oZ\times\oH_{1,4,1},\quad
\oZ\times\oH_{2,4,0},\quad
\oZ\times\oH_{0,4,2}.
\]
(In fact, due to the presence of an extra permutation, this is a single orbit.)
In each of these cases, there are four extra lines
and $\pencil$ extends to a pencil of type $(4,6)$.
Hence, pencils of type $(2,9)$ do not exist, and any pencil of type $(3,7)$
is contained in one of type $(4,6)$.

The value $(p,q)=(1,10)$ can be ruled out similarly, starting with a
subpencil of type $(1,9)$; this time, there are two possibilities
$T=\bH_9(6)$ or $T=\bE_8(6)\oplus\bH_1(6)$.
However, we merely refer to the computation leading to
\autoref{prop.(0,q)}, where pencils of type $(1,10)$ are not excluded
\latin{a priori} but do not appear.
\endproof

\subsection{Proof of \autoref{th.char=3}}\label{proof.char=3}
First, assume that $\Fn X$ is triangle free. If it is also quadrangle free,
we have $\ls|\Fn X|\le44$ by \autoref{lem.quadrangle.free}. If $\Fn X$ has a
line of valency $10$ or more, the bound $\ls|\Fn X|\le37$ is given by
\autoref{prop.(0,q)}. In the remaining case, where $\Fn X$ has a quadrangle
and $\val l\le9$ for each line $l\in\Fn X$, we obtain $\ls|\Fn X|\le49$
directly from~\eqref{eq.quad} and \autoref{lem.quad<=20}.

Now, assume that $\Fn X$ contains a plane $\{a_1,a_2,a_3,a_4\}$.
We use repeatedly the following
special case of~\eqref{eq.pencil.count}, due to
B.~Segre~\cite{Segre}: since
any other line
$l\in\Fn X$
intersects exactly one of $a_i$, $i=1,\ldots,4$ (see \autoref{lem.plane}
and \autoref{rem.rels}),
we have
\[
\ls|\Fn X|=\val a_1+\val a_2+\val a_3+\val a_4-8.
\label{eq.plane}
\]
Configurations containing a pencil of type
$(10,0)$, $(7,0)$, or $(4,6)$ are considered in Corollaries~\ref{cor.(10,0)},
\ref{cor.(7,0)}, and~\ref{cor.(4,6)}, respectively. Otherwise,
by \autoref{cor.exotic.type},
for any
line $l\in\Fn X$ we have $\val l=\ls|\pencil(l)|\le15$.
Hence, $\ls|\Fn X|\le52$ by~\eqref{eq.plane}.
\qed

\remark\label{rem.Fermat}
For an alternative proof of the uniqueness of the Fermat quartic, observe
that, by~\eqref{eq.plane} again, a configuration of size $112$
constitutes a generalized quadrangle $\GQ(3,9)$, see \autoref{s.GQ}, which is
unique up to isomorphism.
As $\rank\bQ_{112}=22$ and
$\discr\bQ_{112}=\<\frac23\>^2$, this lattice admits no
further extension.
\endremark

\section{Pencils in $2$-supersingular quartics}\label{S.2-exotic}
In this section, we discuss
pencils in
supersingular quartics over
an algebraically closed field~$\Bbbk$ of characteristic~$2$
(for short, \emph{$2$-supersingular quartics}).
Our principal observation is the fact that such quartics are related to
indecomposable \emph{odd} negative definite unimodular lattices.

\subsection{The lattice $\NS(X)$}\label{s.NS}
Consider a pencil $\pencil$, not necessarily maximal,
of type $(0,q)$, $q>0$.
First, note that the $3$-torsion of the quotient $\NS(X)/P_{0,q}$ is trivial.
Indeed, by \autoref{lem.discr}, it could be nontrivial only if $q=10$. Then,
$P^\perp$ would be a $2$-elementary lattice,
see \autoref{prop.genus.T}, and, by \autoref{lem.2-discr}, it would
contain a sublattice of the form $\barT(2)$,
where $\barT=\bH_{10}$ or $\bE_8\oplus\bH_2$.
Hence, $P^\perp$, and then also $\NS(X)$, would contain a $(-2)$-vector,
contradicting \autoref{th.Saint-Donat}\iref{Saint-Donat.exceptional}.

Thus, denoting by~$T$ the minimal orthogonal complement and representing
the lattice $\NS(X)$ as an extension of $N:=P_{0,q}\oplus T$ with a certain
kernel~$\CX$ (\cf. \autoref{s.NS.char=3}), from \autoref{prop.genus.T} we
conclude that
the $3$-primary part~$\CX_3$
is the graph of a bijective anti-isometry
\[*
\psi_3\:\discr_3T\longto\discr_3P_{0,q}.
\]
In particular, $\ls|\discr_3T|=9$ and $\discr_3T$ has an
isotropic vector, \viz. the pull-back of~$\Gl$
(see \autoref{obs.discr}); hence, $T$ is an index~$3$
sublattice of a $2$-elementary lattice. By
\autoref{lem.2-discr} again, the latter contains
a sublattice of the form $\barT(2)$, where
$\barT$ is unimodular. Till the rest of this section, we will use the bilinear
form in~$\barT$; in $T$, all products are doubled.

Summarizing,
$\NS(X)$ is a finite index extension of $N:=P_{0,q}\oplus T$, where $T$ can
be described in terms of a unimodular lattice~$\barT$ of rank $20-q$ \via
\[*
T\bigl(\tfrac12\bigr)=\bigl\{a\in\barT\bigm|a\cdot u=0\bmod3\bigr\},\quad
\text{where $u\in\barT$ is fixed, $u^2=q-1\bmod3$}.
\]
The kernel $\CX_3$ is the graph of an anti-isometry
\[
\psi_3\:\discr_3T\longto\discr_3P_{0,q},\qquad
\psi_3\bigl(\barT(2)/T\bigr)=\F_3\Gl.
\label{eq.psi3}
\]
This anti-isometry is unique up to the action of $\OG_h(P_{0,q})$;
hence, with the pair $(\barT,u)$ fixed, the quartic is determined by the
$2$-primary part $\CX_2$.

There is a necessary condition for the
realizability of a pair $(\barT,u)$ by a quartic;
we state it as a separate lemma.

\lemma\label{lem.2-exotic.realizability}
For $(\barT,u)$ as above, the class $u\bmod3\barT$ cannot be represented by
a vector of square $(q-7)$ or a
characteristic vector of square $(q-28)$.
\endlemma

\proof
If $a=u\bmod3\barT$ is as in the statement,
then one of the square~$(-2)$ vectors
$l^*-\sum_{k=1}^q\Gn_k\pm a$ (in the former case) or
$l^*\pm\frac16a$ (in the latter case) is in $\NS(X)$, contradicting
\autoref{th.Saint-Donat}.
\endproof

\proposition\label{prop.2-exotic}
Let $\pencil$ and $(\barT,u)$ be as above, and
assume that the maximal pencil containing~$\pencil$
has exactly $q$ fibers. Then\rom:
\roster
\item\label{2.exotic.odd}
the unimodular lattice~$\barT$ is odd and indecomposable\rom;
in particular, $q\le8$\rom;
\item\label{2.exotic.primitive}
the pencil~$\pencil$ is primitive.
\endroster
\endproposition

\proof
First, we show that $\barT$ does not represent $(-1)$. Indeed, if
$a\in\barT$, $a^2=-1$, then $a\notin T$ by
\autoref{th.Saint-Donat}\iref{Saint-Donat.exceptional};
hence, $a$ represents a
nontrivial class in $\barT(2)/T$ and $\psi_3(a\bmod T)=\pm\Gl$,
see \eqref{eq.psi3}.
It follows that $\Gl\pm a\in\NS(X)$ and, since
we have $(a\pm\Gl)\cdot n_i=0$
for $i=1,\ldots,q$, the
pencil has an extra fiber.

If $q=12$ and $\barT=\bE_8$,
we can assume that $u^2=-4$; then
the vector $-2a$ is characteristic,
contradicting \autoref{lem.2-exotic.realizability}.
With $\bE_8$ eliminated,
we have $q\le8$ and $\barT$ is automatically odd unless $q=4$. However,
in the latter case, $\pencil$ is primitive, see
\autoref{lem.discr}, and $\discr_2T$ must be odd; hence, $\barT$ is also odd.
Now, it follows that $T$ is indecomposable, as the smallest decomposable odd
unimodular lattice not representing~$(-1)$ is $\bE_8\oplus\bD_{12}^+$ of
rank~$20$ (corresponding to $q=0$).

Similarly, by \autoref{lem.discr}, $\pencil$ is primitive unless $q=8$, in
which case the only admissible extension~$\tP$ has even discriminant
$\discr_2\tP$. Then $\NS(X)$ contains $\frac12a$, where
$a\in\barT=\bD_{12}^+$ is a characteristic vector.
One can check that, for any choice of $u\bmod3\barT$,
there is a characteristic vector
$a\in u^\perp$ of square~$(-4)$ (see \autoref{ss.q=8} below for a detailed
description of this lattice).
\endproof

\subsection{Orbits}\label{s.orbits}
Due to \autoref{prop.2-exotic}\iref{2.exotic.primitive}, the kernel
$\CX_2\subset\discr_2N$
is  the graph of an anti-isometry
\[*
\psi\:\CH\longto\discr_2P_{0,q},\qquad \CH\subset\discr_2T=\barT/2\barT.
\]
Denote by $c_P\in\discr_2P_{0,q}$ and $c_T\in\discr_2T$ the characteristic
vectors; they are both nontrivial. Since $\discr_2\NS(X)$ must be even, we
have
\[*
c_T\in\CH=\operatorname{Domain}\psi\quad\text{and}\quad
\psi(c_T)=c_P.
\]

Consider the group
\[*
\OG(\barT,u):=\bigl\{g\in\OG(\barT)\bigm|g(u)=\pm u\bmod3\barT\bigr\}.
\]
As in the case of characteristic~$3$, the orbits of the $\OG_h(N)$-action on
$\discr_2N$ split into products
\[*
\oN_k\times\ooO,\quad
 \text{where $k=0,\ldots,q$ and $\ooO\subset\barT/2\barT$ is an
 $\OG(\barT,u)$-orbit}.
\]
Furthermore, we only need to check
the admissibility and compute the number of extra lines, \ie, the fibers of
the natural maps
\[*
\sec_2\:\{\text{sections of~$\pencil$}\}\longto\CX_2,\qquad
\fib_2\:\pencil(l)\sminus\pencil\longto\CX_2,
\]
for one
representative of each orbit.
Then, an isotropic subgroup $\CX_2\subset\discr_2N$
is admissible if and only if so are all its elements, and the number of lines is
additive.
Note that $N_q=\{c_P\}$ and
the orbit $\oN_q\times\{c_T\}$ is always in~$\CX_2$.

\remark
In the case of characteristic~$2$, we have a better control over the geometry
of the extra lines.
Denote by
\[*
\pr\:\discr_2N\longto\discr_2P_{0,q}
\]
the projection.
Then, for an extra line $s\in\NS(X)\sminus N$, we have
\[*
\pr(s\bmod N)=\sum_{k\in I}3\Gn_k,\quad
\text{where $I:=\{k\in\fb_1(\pencil)\,|\,s\cdot n_k=1\}$}.
\]
As a consequence, all nonempty fibers of the composition ${\pr}\circ\fib_2$ are
over~$\oN_1$:
a generator $3\Gn_k$ is in the image if and only if the $k$-th
$1$-fiber of~$\pencil$ becomes a $3$-fiber. Alternatively, all extra lines in
$\pencil(l)\sminus\pencil$
are of the form
\[*
n_k+3\Gn_k\pm\tfrac12a=-\Gl+\Gn_k\pm\tfrac12a,
\]
$a\in\barT$,
$a^2=-3$, provided that these vectors are in $\NS(X)$, \cf.
\autoref{lem.psi2}.
\endremark

\subsection{The values $q=5$, $6$, and~$8$}\label{s.large-q}
In this section, we describe the isomorphism classes of pairs $(\barT,u)$
appearing in \autoref{prop.2-exotic} for the large values $q=5$, $6$,
and~$8$. (Note that $q\ne7$, as there is no indecomposable lattice of
rank~$13$.) There are eight classes, listed below.
Afterwards, it is straightforward,
although tedious, to use \GAP~\cite{GAP4} and enumerate all
kernels~$\CX_2$; we merely state the final result in
\autoref{prop.2-exotic.counts}.
An important experimental fact is that $T$ is necessarily
primitive in~$\NS(X)$; hence, instead of
$\psi\:\CH\to\DD\subset\discr_2P_{0,q}$, we can consider its
inverse $\psi\1\:\DD\to\discr_2T$, which is defined on a subspace
$\DD\subset{-\CS_q}$ (see
\autoref{s.discr.lemmas}) containing the characteristic vector.
Up to the action of~$\SG{q}$ by
permutations of the generators, which is induced from $\OG(P_{0,q})$,
there are relatively few such subspaces.

\subsubsection{The case $q=8$ and $\barT=\bD_{12}^+$}\label{ss.q=8}
We represent $\bD_{12}$ as the maximal even sublattice in $\bH_{12}$.
This lattice has three unimodular extensions, all odd: one is
the original lattice $\bH_{12}$,
and the two others are isomorphic. We denote by
$\bD_{12}^+$ the extension by $\frac12\be$; then,
$\OG(\bD_{12}^+)\subset\OG(\bD_{12})$ is the
index~$2$ subgroup generated by
reflections against vectors of square~$(-2)$.
Since $u^2=1\bmod3$, the $\OG(\bD_{12}^+)$-orbits of
vectors $u\bmod3\bD_{12}^+$
are characterized by the Hamming norm: we have $\oH_2$, $\oH_5$,
$\oH_8$, and~$\oH_{11}$.
(Note, in particular, that there always is a characteristic
vector, \viz. $2\e_{12}\in u^\perp$,
of square~$(-4)$; this fact was used in the proof of
\autoref{prop.2-exotic}.) The only orbit satisfying
\autoref{lem.2-exotic.realizability} is~$\oH_8$.

\subsubsection{The case $q=6$ and $\barT=(\bE_{7}^2)^+$}\label{ss.q=6}
This lattice is the only nontrivial extension of $\bE_7^2$; its kernel is
generated by the vector of square~$1$ in $\discr\bE_7^2=\<\frac12\>^2$.
The orbits of the action of $\OG(\bE_7)$ on $\bE_7/3\bE_7$ are almost
distinguished by the length of the shortest representatives:
we have
\[*
\ooE_r:=\bigr\{u\bmod3\bE_7\bigm|u\in\bE_7,\ u^2=-r\bigr\},\quad
 r=0,2,4,6,8,
\]
but the set $\ooE_6$ splits into two orbits $\ooE_6'$, $\ooE_6''$,
where $\ooE_6''$ is
characterized by the fact that its shortest representatives vanish
$\bmod\,2\bE_7\dual$. Since $u^2=2\bmod3$ (and in view of the obvious
symmetry), we have the $\OG(\barT)$-orbits
\roster*
\item
$\ooE_2\times\ooE_8$, $\ooE_4\times\ooE_0$, $\ooE_4\times\ooE_6''$, which
contradict \autoref{lem.2-exotic.realizability}, and
\item
$\ooE_2\times\ooE_2$, $\ooE_8\times\ooE_8$, $\ooE_4\times\ooE_6'$, which
give rise to nonsingular quartics.
\endroster

\subsubsection{The case $q=5$ and $\barT=\bA_{15}^+$}\label{ss.q=5}
This is the only unimodular extension of $\bA_{15}$, \viz. the one by
$4\discr\bA_{15}$.
Represent~$\bA_{15}$ as $\be^\perp\subset\bH_{16}$.
Then, $\OG(\bA_{15}^+)=\OG(\bA_{15})$ is the subgroup of $\OG(\bH_{16})$
stabilizing the set $\{\pm\be\}$.
The orbits of the action of this group
on the $\F_3$-vector space $\bA_{15}/3\bA_{15}$ are
\[*
\oH_{r,s}\ni 2(\e_1+\ldots+\e_r)+(\e_{r+1}+\ldots+\e_{r+s})
 -(\e_{r+s+1}+\ldots+\e_{3r+2s});
\]
we have $r,s\ge0$, $3r+2s\le16$, and $s=1\bmod3$ (since $u^2=1\bmod3$).
Five orbits contradict \autoref{lem.2-exotic.realizability};
the remaining valid orbits are
$\oH_{0,4}$, $\oH_{0,7}$, $\oH_{1,4}$, $\oH_{3,1}$.

\proposition\label{prop.2-exotic.counts}
Let $X\subset\Cp3$ be a $2$-supersingular quartic, and let
$\pencil\subset\Fn X$ be a maximal pencil of type $(p,q)$.
Then the number $p+q$ of fibers of~$\pencil$
can take only one of the following
values\rom:
\roster*
\item
$p+q=8$\rom: then $p\le3$ and $\ls|\Fn X|=40$ or $\ls|\Fn X|\le32$\rom;
\item
$p+q=6$\rom: then $p\le5$ and $\ls|\Fn X|\le32$\rom;
\item
$p+q=5$\rom: then $\ls|\Fn X|\le24$\rom;
\item
$p+q\le4$.
\endroster
If $p+q=8$,
there are
{\3at most}
four distinct configurations
$\Fn X$ of size~$40$\rom; they are as in
\autoref{th.char=2}\iref{40.1}--\iref{40.5}.
\endproposition

{\3Hypothetically,
the bounds in \autoref{prop.2-exotic.counts} are sharp:
they are realizable by homological configurations, \cf.
\autoref{rem.Saint-Donat}.}
In the spirit of Propositions~\ref{prop.types} and~\ref{cor.exotic.type},
the part of \autoref{prop.2-exotic.counts} concerning the types of the
pencils can be summarized in the form of the table
\[*
\centerline{\setbox0\hbox{00}\vbox{\halign{
 \hss\strut$#$&&\quad\hbox to\wd0{\hss$#$\hss}\cr
p=  &5&4&3&2&1&0\cr
q\le&1&2&5&6&7&8\rlap,\cr}}}
\]
with the extra restriction that $p+q\ne7$.

\subsection{The value $q=4$}\label{s.q=4}
The only indecomposable odd unimodular lattice of rank~$16$ is
$\barT=(\bD_8^2)^+$.
Recall that $\bD_8$ is the maximal even sublattice in $\bH_8$. Then,
$(\bD_8^2)^+$ can be described as the extension of
$\bD_8^2$ by the two $({\bmod}\,\Z)$-orthogonal vectors $\frac12\be\oplus \e_8$,
$\e_8\oplus\frac12\be$ of square $1\bmod2\Z$.
The $\OG(\bD_8)$-orbits on $\bD_8/3\bD_8$ are almost
characterized by the Hamming norm: we have $\oH_0$ through $\oH_7$, and
$\oH_8$ splits into
two orbits $\oH_8^+\ni\be$, $\oH_8^-\ni\be-2\e_8$.
Hence, with the obvious symmetry taken into account, there are $17$ orbits of
vectors $u\bmod3\barT$ of square $0\bmod3$. Nine of them contradict
\autoref{lem.2-exotic.realizability}, and the remaining eight are
\roster*
\item
$\oH_4\times\oH_5$, $\oH_4\times\oH_8^+$, $\oH_7\times\oH_5$,
$\oH_3\times\oH_6$, $\oH_6\times\oH_6$, and
\item
$\oH_4\times\oH_2$, $\oH_7\times\oH_8^-$, $\oH_0\times\oH_3$.
\endroster
The last three orbits are characterized by the fact that a generic quartic
has at least one section; this section intersects all four fibers
(in particular, the configuration of lines is not quadrangle free).

Unlike the three cases considered in the previous section, this time we have
$\dim(\CX_2\cap\discr_2T)\le3$.
Together with the size of the spaces and groups involved, this fact
makes the computation difficult. For this reason, we only consider two extremal
cases: quadrangle free configurations and those where the maximal pencil
containing~$\pencil$ has at least three $3$-fibers.

\proposition\label{prop.2-exotic.q=4}
Let $X\subset\Cp3$ be a $2$-supersingular quartic, and let
$\pencil\subset\Fn X$ be a maximal pencil of type $(p,q)$, $p+q=4$.
Then\rom:
\roster*
\item
if $p\ge3$, then $\ls|\Fn X|=40$ or $\ls|\Fn X|\le32$\rom;
\item
if $\Fn X$ is quadrangle free, then $\ls|\Fn X|\le9$.
\endroster
\endproposition

\remark\label{rem.char=2.counts}
The {\3potential} line counts $\ls|\Fn X|$
observed in the course of the proof of
Propositions~\ref{prop.2-exotic.counts} and~\ref{prop.2-exotic.q=4} are
\[*
\bigl\{5,6,\ldots,17, 18, 20, 22, 24, 28, 32, 40\bigr\}.
\]
This list
is complete for the
{\3homological (\cf. \autoref{rem.Saint-Donat})}
configurations containing a pencil with at least five
fibers.
In more detail, for a nonsingular quartic $X\subset\Cp3$
we have $\Artin(X)\ge3$ and the
values of $\ls|\Fn X|$ are distributed as follows:
\roster*
\item
if $\Artin(X)=3$, then $\ls|\Fn X|=40$ (five quartics)\rom;
\item
if $\Artin(X)=4$, then $\ls|\Fn X|\in\{12,16,20,24,28,32\}$\rom;
\item
if $\Artin(X)=5$, then $|\Fn X|\le24$ and $\ls|\Fn X|=0\bmod2$\rom;
\item
if $\Artin(X)=6$, then $|\Fn X|\le18$ or $\ls|\Fn X|=20$\rom;
\item
if $\Artin(X)=7,8$, then $|\Fn X|\le16,12$, respectively.
\endroster
This observation substantiates and refines \autoref{conj.char=2}.
\endremark


\subsection{Proof of \autoref{th.char=2}}\label{proof.char=2}
In view of \autoref{prop.2-exotic.counts}, we can
assume that each pencil $\pencil\subset\Fn X$ has at most four fibers; in
particular, $\ls|\pencil|\le12$.

First, assume that $\Fn X$ is triangle free. If $\Fn X$ is also quadrangle
free, then, by \autoref{prop.2-exotic.q=4}, either $\Fn X\le32$ or $\Fn X$ is
locally elliptic; in the latter case, $\Fn X\le31$ by
\autoref{lem.locally.elliptic}. If $\Fn X$ has a quadrangle, then each
valency $\val l_i\le4$
in~\eqref{eq.quad} and, even if $\flines{\pi_Q}=24$,
we have $\ls|\Fn X|\le32$ again.

Now, assume that $\Fn X$ has a plane.
Configurations containing a pencil of type $(p,q)$ with
either $p+q\ge5$ or $p+q=4$ and $p\ge3$ are considered, respectively, in
Propositions~\ref{prop.2-exotic.counts} and~\ref{prop.2-exotic.q=4}.
Otherwise, $\val l=\ls|\pencil(l)|\le9$ for any
line $l\in\Fn X$ and, hence, $\ls|\Fn X|\le28$ by~\eqref{eq.plane}.
(In fact, $\ls|\Fn X|\le27$, as there are no generalized quadrangles
$\GQ(3,2)$, see \autoref{s.GQ}.)

There remains to show that each configuration $\Fn X$ of size~$40$
is one of those
listed in the statement.
In view of \autoref{prop.2-exotic.counts}, we can
assume that each pencil $\pencil\subset\Fn X$ has at most four fibers; in
particular, $\ls|\pencil|\le12$. Then, by~\eqref{eq.plane} again,
$\Fn X$
constitutes
one of the two generalized quadrangles $\GQ(3,3)$, see \autoref{s.GQ}, which
we consider separately.
Let~$T$ be the orthogonal complement of~$\bQ_{40}^*$ in $\NS(X)$: it is a
negative definite lattice of rank~$6$.

Let $\Fn X\cong\Q(4,3)$.
Since
$\discr\bQ_{40}'=\<\frac43\>\oplus\CV_2^3$,
we have $\discr T=\<\frac23\>\oplus\CU_2^r$, $r\le3$, see \autoref{th.NS} and
\autoref{prop.genus.T}.
In fact, $r\le2$ by \autoref{th.Nikulin}, and, by an analogue of
\autoref{lem.2-discr}, $T$ is an extension of $\bD_4\oplus\bA_2(2)$,
which is unique in its genus.
This lattice has $(-2)$-vectors, contradicting
\autoref{th.Saint-Donat}\iref{Saint-Donat.exceptional}.
Alternatively, one can start with the lattice~$T'$ given by
\autoref{lem.2-discr}.
Then, $\discr_2T'$ is necessarily odd;
hence, $T'=\bA_2(2)\oplus\bA_1^4$ and $T\supset T'$ is the extension by the
characteristic vector $c\in\discr_2T'$.
This observation proves both the uniqueness in the genus and the existence of
$(-2)$-vectors, as such vectors are already present in~$T'$.

Let $\Fn X\cong W(3)$.
Since $\discr\bQ_{40}''=\<\frac23\>^5$, we have
$\discr T=\<\frac43\>^5\oplus\CU_2^r\oplus\CV_2$ for some $r\le2$.
The lattice $\bQ_{40}''$ has no admissible
finite index extensions and, by an analogue of
\autoref{lem.2-discr}, $T$ is an extension of $\bbE_6(2)$, which is unique in
its genus.
Since the natural homomorphism $\OG(\bbE_6)\to\Aut\discr_3\bbE_6$
is an isomorphism,
an anti-isometry $\discr_3\bQ_{40}''\to\discr_3\bbE_6(2)$ is
essentially unique, and
we do obtain a quartic~$X$ as in \autoref{th.char=2}\iref{40.5}. A simple
computation shows that the lattice $\NS(X)$ has no admissible
finite index extensions.
\qed

\subsection{Other generalized quadrangles}\label{s.GQ.X}
For completeness,
we discuss briefly the realizability (in the sense of \autoref{s.GQ}) of the other
generalized quadrangles
$\GQ(3,t)$ by configurations of lines in nonsingular quartics.

\subsubsection{The quadrangle $\GQ(3,9)$}\label{ss.GQ9}
The only realization of $\GQ(3,9)\cong Q(5,3)$ is that by the Fermat quartic in
characteristic~$3$, see \autoref{th.char=3}\iref{112.(10,0)} and
\autoref{rem.Fermat},
since for all other quartics~$X$ one has $\ls|\Fn X|\le64$.
In the Fermat quartic~$X$, the four lines constituting a plane
$\Ga$ intersect at a single point $P_\Ga$, see
\autoref{cor.4-point} below. Hence, taking $\Fn X$ for the set of lines and
all points $P_\Ga$ for the set of points, we also obtain a generalized
quadrangle, one of type $\GQ(9,3)$.

\subsubsection{The quadrangle $\GQ(3,5)$}\label{ss.GQ5}
This generalized quadrangle is not realized by a configuration of lines.
We have
$\rank\bQ_{64}=19$ and $\discr\bQ_{64}=\CU_2^2\oplus\<\frac74\>$, and
this lattice has no admissible finite index extensions.
Thus, according to \autoref{th.realization}, this configuration could only be realized by
a $2$-supersingular quartic, which would contradict
\autoref{prop.2-exotic.counts}.

\subsubsection{The quadrangles $\GQ(3,3)$}\label{ss.GQ3}
As shown in \autoref{proof.char=2}, only
the generalized quadrangle $W(3)$ appears as the configuration
of lines in a supersingular quartic~$X$ over a field of characteristic~$2$.
If $X$ is
supersingular in characteristic~$3$, both $Q(4,3)$ and $W(3)$ can be
realized:
\roster
\item\label{S(4,3)}
if $\Fn X\cong Q(4,3)$, then $\Artin X=3$ and $(\bQ_{40}')^\perp=\bbE_6(2)$;
\item\label{W(3)}
if $\Fn X\cong W(3)$, then $3\le\Artin X\le5$ and $(\bQ_{40}'')^\perp=\bbE_6$.
\endroster
In case~\iref{S(4,3)}, the representation is unique, as any proper admissible
extension of the lattice $\NS(X)=\bQ_{40}'\oplus\bbE_6(2)$
contains extra lines. In case~\iref{W(3)},
the configuration is maximal with respect to inclusion:
the minimal lattice $\NS(X)=\bQ_{40}''\oplus\bbE_6$ has admissible extensions
of index~$3$ and~$9$, all
with the same configuration of lines.

Both $Q(4,3)$ and $W(3)$ are realizable by nonsingular quartics over~$\C$,
with the transcendental
lattice isomorphic to $\bA_2(2)\oplus\bU^2(2)$ and $\bA_2\oplus\bU^2(3)$,
respectively. The former contains $\bU(2)$, and the latter contains
$\bA_1$; hence, both quadrangles can be represented by \emph{real}
lines in a real quartic, see~\cite[Corollary 3.14]{DIS}.

\subsubsection{The quadrangle $\GQ(3,1)$}\label{ss.GQ1}
We have $\rank\bQ_{16}=10$ and $\discr\bQ_{16}=\CU_4^2$.
Hence,
the quadrangle
$\GQ(3,1)\cong Q(3,3)$ can be realized
by a nonsingular quartic over~$\C$, with
the transcendental lattice $\bE_8\oplus\bU^2(4)\supset\bA_1$;
by~\cite{DIS}, both the quartic and the
lines can be chosen real. We omit the discussion of the realizability
of
$\GQ(3,1)$
by supersingular quartics in characteristics~$2$ and~$3$.



\section{Geometric arguments}\label{S.geometry}

In this section, we employ direct geometric arguments (at the level of
defining equations) to establish the uniqueness of Schur's quartic and prove
\autoref{th.ordinary}.

\subsection{Pairs of obverse pencils}\label{s.pair}
Let $\pencil_i:=\pencil(l_i)$, $i=1,2$, $l_i\in\Fn X$, be a pair of obverse
pencils. The pencil $\pencil_i$ defines a (quasi-)elliptic pencil
$\pi_i\:X\to\Cp1$; hence, we have a map
$\pi[l_1,l_2]:=\pi_1\times\pi_2\:X\to\Cp1\times\Cp1$. The base $\Cp1$ of the
projection~$\pi_i$ can be identified with~$l_j$, $j\ne i$. Hence, the
pull-back of a point $(x,y)\in\Cp1\times\Cp1$ is the intersection of~$X$
with the line connecting the points $y\in l_1$ and $x\in l_2$, excluding
$x,y$ themselves. It follows that $\pi:=\pi[l_1,l_2]$ is of degree $2$.
The deck translation of the double covering $X\to\Cp1\times\Cp1$ is known as
the \emph{Segre involution}; typically, it is not projective.

First, assume that $\kchar\ne2$.
Then $\pi$ is a double covering ramified at a
curve $D\subset\Cp1\times\Cp1$ of bidegree~$(4,4)$.
We keep the notation $(x,y)$ for the coordinates in
the target quadric $\Cp1\times\Cp1$.

\observation\label{obs.covering}
The following statements are straightforward
(where 
a curve is called \emph{even}
if its
intersection index with~$D$ at each intersection point is even):
\roster
\item\label{double.l1}
the line $l_1$ projects to an even irreducible curve of bi-degree $(3,1)$;
\item\label{double.l2}
the line $l_2$ projects to an even irreducible curve of bi-degree $(1,3)$;
\item\label{double.common}
a line $a\in\pencil_1\cap\pencil_2$ contracts to a singular point of~$D$;
\item\label{double.nodes}
any singular point of~$D$ is
a simple node of the form $\pi(a)$,
$a\in\pencil_1\cap\pencil_2$;
\item\label{double.sing}
the curves $\pi(l_1)$ and $\pi(l_2)$ contain
all points $\pi(a)$, $a\in\pencil_1\cap\pencil_2$;
\item\label{double.y}
a line $a\in\pencil_1\sminus\pencil_2$ projects to an even generatrix
$y=\const$;
\item\label{double.x}
a line $a\in\pencil_2\sminus\pencil_1$ projects to an even generatrix
$x=\const$;
\item\label{double.other}
any other line $b$ projects to an even irreducible curve of
bidegree $(1,1)$; this curve contains a point $\pi(a)$,
$a\in\pencil_1\cap\pencil_2$,
if and only if $b$ intersects~$a$.
\endroster
\endobservation

As an immediate consequence, the projection~$\pi$ establishes a
canonical one-to-one correspondence
between the $3$-fibers of~$\pencil_1$
(respectively, $3$-fibers of~$\pencil_2$) and the even generatrices of the form
$y=\const$ (respectively, $x=\const$) passing through a singular point
of the ramification locus~$D$.

\proposition\label{prop.10-2}
Assume that $\kchar\ne2$ and that $\Fn X$ contains a configuration
$a_1,\ldots,a_{10},b_1,b_2$ as in \autoref{lem.10-2}. Then, the ramification
locus~$D$ of
$\pi:=\pi[b_1,b_2]$
splits into the union $\pi(b_1)\cup\pi(b_2)$.
\endproposition

\proof
According to \autoref{obs.covering}\iref{double.common}, \iref{double.sing},
the ramification locus~$D$ must have ten simple nodes, which must coincide
with the ten points of $\pi(b_1)\cap\pi(b_2)$.
Now, the statement is an immediate consequence of B\'{e}zout's theorem and
the fact that both $\pi(b_1)$ and $\pi(b_2)$ are irreducible.
\endproof

\corollary\label{cor.first.kind}
Under
the assumptions of \autoref{prop.10-2}, the lines $b_1$, $b_2$ are of
the first kind in the sense of Segre~\cite{Segre}.
\endcorollary

\proof
Recall that a line~$l$ is said to be of the second kind if a generic
fiber~$C$ of the pencil $\pi[l]$ is a nonsingular cubic and the points of
intersection of~$l$ and~$C$ are inflection points of~$C$. Assume that
the line~$b_1$
is of the second kind and let $b_1\cap C=\{P_1,P_2,P_3\}$.
Then the three differences $P_i-P_j$, $1\le i<j\le3$, are $3$-torsion points
of the Jacobian $J(C)$. On the other hand, since $b_1$ is contained in the
ramification locus of the fiberwise degree~$2$ map $\pi[b_1,b_2]$, these
differences are $2$-torsion points of $J(C)$. This is a contradiction.

Alternatively, one can
show that, typically, $P_1$, $P_2$, $P_3$ are not inflection points
by computing the Hessian of a generic fiber
in~\eqref{eq.10-2} below.
\endproof

Under the assumptions of \autoref{prop.10-2}, the ramification locus
$\pi(b_1)\cup\pi(b_2)$ determines both the abstract $K3$-surface~$X$ and
polarization~$h$, \eg, as the sum of the reduced preimage of $\pi(b_1)$ (the
line~$b_1$) and the pull-back of any fiber $y=\const$
(see \autoref{lem.plane}).
Hence, it determines
the quartic $X\subset\Cp3$, and one can easily see that the latter is given
by the equation
\[
z_1^3z_3f_1\Bigl(\frac{z_0}{z_1},\frac{z_2}{z_3}\Bigr)
=z_1z_3^3f_2\Bigl(\frac{z_0}{z_1},\frac{z_2}{z_3}\Bigr),
\label{eq.10-2}
\]
where $f_i(x,y)=0$ is a defining equation of the component $\pi(b_i)$,
$i=1,2$.
Clearly, this equation can be rewritten in the form
\[*
\bar f_1(z_0,z_1;z_2,z_3)=\bar f_2(z_0,z_1;z_2,z_3),
\]
where $\bar f_i(z_0,z_1;z_2,z_3)$ is the homogenization of $f_i(x,y)$,
$i=1,2$.
Conversely, given two polynomials $f_1$, $f_2$ of bidegree $(3,1)$ and
$(1,3)$, respectively, then,
assuming that $\kchar\ne2$, the quartic
given by~\eqref{eq.10-2} is nonsingular if
and only if the two curves $B_i:=\{f_i(x,y)=0\}\subset\Cp1\times\Cp1$ are
nonsingular (equivalently, irreducible) and intersect transversally at ten
points.

\proposition\label{prop.10-2.lines}
Under the assumptions of \autoref{prop.10-2}, the lines $l\in\Fn X$ disjoint
from $b_1,b_2$ are in a two-to-one correspondence with quadruples
of points $A_i\in\pi(b_1)\cap\pi(b_2)$, $i=1,\ldots,4$, whose coordinates
$(x_i,y_i)$ have equal unharmonic ratios\rom:
$(x_1,x_2;x_3,x_4)=(y_1,y_2;y_3,y_4)$.
\endproposition

\proof
By \autoref{lem.10-2} (see also
\autoref{rem.rels}), any line $l\in\Fn X$ disjoint from $b_1,b_2$ intersects
exactly four of the ten lines $a_1,\ldots,a_{10}$; hence, its image is
an irreducible bidegree $(1,1)$ curve~$L$ passing through the four
images $A_i:=\pi(a_i)\in\pi(b_1)\cap\pi(b_2)$.
Conversely, any such curve
is even and each of the two components of its pull-back intersects
a generic plane at a single point.
\endproof

Now, assume that $\kchar=2$.
Then, in appropriate coordinates in $\Cal O_{\Cp1\times\Cp1}(4,4)$,
the surface~$X$ is given by
\[*
z^2+f_2(x,y)z+f_4(x,y)=0,
\]
where $f_2$ and $f_4$ are of bidegree $(2,2)$
and $(4,4)$, respectively.
One has $f_2\ne0$, \ie, the covering cannot be purely
inseparable.
Indeed,
otherwise, the two pencils $\pi[l_1]$, $\pi[l_2]$ would be
quasi-elliptic, which would contradict \autoref{prop.Euler}.
Thus, the formal ramification locus of~$\pi$ is the non-reduced
curve $D=2\bar D$, where $\bar D$ is given by $f_2=0$ and has
bidegree~$(2,2)$. Hence, any point of~$D$ is singular and any curve in
$\Cp1\times\Cp1$ is even. With this understood, all statements
except~\iref{double.nodes} in \autoref{obs.covering} hold literally, although
most are no longer invertible.

\proposition\label{prop.10-2.char2}
If $\kchar=2$, then $\Fn X$ cannot contain a configuration as in
\autoref{lem.10-2}, \ie, two
disjoint
lines~$l_1$, $l_2$ can intersect at most eight other lines.
\endproposition

\proof
The projection $\pi:=\pi[l_1,l_2]$ contracts any line
$a\in\pencil(l_1)\cap\pencil(l_2)$
to a point $A\in\Cp1\times\Cp1$ common to the
bidegree~$(2,2)$ curve~$\bar D$ and irreducible bidegree~$(3,1)$ curve
$\pi(l_1)$. By B\'{e}zout's theorem, there are at most eight such points.
\endproof

\subsection{Schur's quartic}\label{s.Schur}
According to~\cite{DIS}, if $\kchar=0$, there is a unique quartic $X_{64}$
containing $64$ lines; it is called \emph{Schur's quartic}. The lines
contained in~$\X_{64}$ span the N\'{e}ron--Severi lattice $\NS(X)$; we denote
this lattice by $\bX_{64}$ and call it the \emph{Schur configuration}.
The following property of~$\bX_{64}$ is easily
deduced from the explicit incidence matrix of $\Fn\bX_{64}$ or,
alternatively, from the classical description of the $64$ lines in Schur's
quartic, see, \eg,~\cite{Barth:K3}

\lemma\label{lem.Schur.lines}
The configuration $\Fn\bX_{64}$ contains a collection
$a_1,\ldots,a_{10},b_1,b_2,c$ with the following properties\rom:
\roster
\item\label{Schur.10-2}
the lines $a_1,\ldots,a_{10},b_1,b_2$ are as in \autoref{lem.10-2}\rom;
\item\label{Schur.3-fibers}
the four lines $a_1,\ldots,a_4$ are in the $3$-fibers of both~$\pencil(b_1)$
and~$\pencil(b_2)$\rom;
\item\label{Schur.c}
the line~$c$ intersects $a_1,\ldots,a_4$.
\endroster
\endlemma

\lemma\label{lem.Schur.uniqueness}
A nonsingular quartic $X\in\Cp3$
whose configuration of lines $\Fn X$ has the properties stated
in \autoref{lem.Schur.lines} is unique up to projective transformation. Such
a quartic exists if and only if $\kchar\ne2$.
\endlemma

\proof
The possibility $\kchar=2$ is ruled out by \autoref{prop.10-2.char2}.
Henceforth, we assume that $\kchar\ne2$, and the existence of~$X$ is
given by an explicit example, \viz. Schur's quartic.
Below, without further references, we freely use
\autoref{obs.covering}.

Let $\pi:=\pi[b_1,b_2]$.
According to \autoref{prop.10-2}, the
ramification locus~$D$ of~$\pi$ splits into two irreducible curves
$B_i:=\pi(b_i)$, $i=1,2$, which intersect at the ten points $A_i:=\pi(a_i)$,
$ i=1,\ldots,10$. Furthermore, the irreducible bidegree $(1,1)$ curve
$C:=\pi(c)$ passes through the four points $A_1,\ldots,A_4$; it follows that
the coordinates of these points can be chosen to be $(r,r)$, where
$r=0$, $1$, $\infty$, or $\Gl\in\Bbbk\sminus\{0,1\}$.

An additional property of the irreducible bidegree $(3,1)$ curve $B_1$ is
that it is tangent to the four lines $y=r$ (as a special case, inflection
tangent at $x=r$). To find this curve, we change the $x$-coordinate so that
the first three \emph{tangency} points are at $x=0,1,\infty$, respectively.
Then, the equation has the form
\[*
x^3+(u-1)x^2+y(v_1x+v_0)=0,
\]
and the tangency at $(1,1)$ gives us $v_0 = u+1$, $v_1 = -2u-1$.

First, assume that the degree~$3$ projection $B_1\to\Cp1$ is not purely
inseparable.
Then, the
fourth tangency level is
\[*
y=\Gl(u)=-\frac{(u-1)^3(u+1)}{(2u+1)^3};
\]
since $\Gl\ne0,1,\infty$, we have $u\ne0,\pm1,-2,-\frac12$. In these
coordinates, the additional points of intersection of~$B_1$ with the
four tangents $y=r$ are
\[*
x_0=-u+1,\quad x_1=-u-1,\quad
x_\infty=\frac{u+1}{2u+1},\quad
x_\Gl=\frac{u-1}{2u+1}.
\]
Now, switching back to the original coordinates, \ie, sending $x_0,x_1,x_\infty$
to $0,1,\infty$, respectively, we arrive at the following $x$-coordinate
of~$A_4$:
\[*
\mu(u):=-\frac{(u+1)^3(u-1)}{2u+1}.
\]
Equating $\mu(u)=\Gl(u)$ and disregarding the values of~$u$ ruled out above,
we obtain $u^2+u+1=0$,
\ie, $u=\epsilon_{1,2}$ is a primitive $3$-rd root of unity.
In particular, since $u\ne1$, we have $\kchar\ne3$.

If $B_1\to\Cp1$ is purely inseparable, then $\kchar=3$ and $u=1$, \ie,
$B_1$ is given by $x^3=y$.
In this case, one easily
finds that $\Gl=\Gm=-1$.

The curve~$B_2$ has similar properties,
and it
passes through the same quadruple of points $A_1,\ldots,A_4$; hence,
its equation is that of $B_1$, {\em with the same value} of $u=1$
(if $\kchar=3$),
$\epsilon_1$, or~$\epsilon_2$, with $x$ and $y$ interchanged.
It is easily seen that the two
values $u=\epsilon_{1,2}$ can be interchanged by an appropriate change of coordinates.
Thus, the ramification locus $D\subset\Cp1\times\Cp1$ is unique up to
isomorphism, and
the quartic $X\in\Cp3$ is given by~\eqref{eq.10-2} with $f_1,f_2$ as
described above.
\endproof

\corollary\label{cor.Schur}
A nonsingular quartic $X\subset\Cp3$
with $\Fn X\cong\Fn\bX_{64}$ exists if and only if
$\kchar\ne2$ or~$3$. If exists, $X$ is isomorphic to Schur's quartic.
\endcorollary

\proof
The uniqueness and the restriction $\kchar\ne2$ are given by
\autoref{lem.Schur.uniqueness}. If $\kchar\ne2$ or~$3$, the classical
Schur quartic is nonsingular and contains exactly $64$ lines, which follows
from the classical description of these lines. If $\kchar=3$,
we obtain $\epsilon=1$ in the proof of \autoref{lem.Schur.uniqueness} and
both pencils become quasi-elliptic; hence, $X$ becomes supersingular and,
by \autoref{th.char=3}, $X$ is the Fermat quartic with $112$ lines.
The same
conclusion can as well be derived from the explicit equation.
\endproof

It is fairly easy to describe all quartics satisfying
conditions~\iref{Schur.10-2}, \iref{Schur.3-fibers} of
\autoref{lem.Schur.lines}, \ie, without assuming the existence of line~$c$.
To this end,
we should not assume that $A_4$ has equal coordinates
in the proof of \autoref{lem.Schur.uniqueness}.
Thus, we merely start with a pair of curves~$B_1$, $B_2$
with distinct values~$u$, $v$ of the parameter and
equate $\Gl(u)=\Gm(v)$ and
$\Gl(v)=\Gm(u)$. This gives us a $1$-parameter family
\[*
2uv+u+v+2=0
\]
and a number of discrete pairs, satisfying
\[*
(2v+1)^2u^3+3v(2v+1)u^2-3v(2v+1)u-v(v+2)^2=0
\]
and the same equation with $u$, $v$ interchanged. The two latter result in
three sets of Galois conjugate solutions:
\roster*
\item
$2u^4+4u^3+2u^2+1=0$ and $v=-u-1$,
\item
$u^4+2u^2+4u+2=0$ and $v=-u^3+u^2-3u-2$, and
\item
$u^4+4u^3+8u^2+4u+1=0$ and $v=-u^3-4u^2-8u-4$,
\endroster
probably not all distinct. Strictly speaking, with $u$ fixed, over some
primes there may be other solutions for~$v$.
However, we do not investigate this issue
any further, nor do we discuss the conditions under which the ramification
loci obtained do give rise to nonsingular quartics.

\subsection{Proof of \autoref{th.ordinary}}\label{proof.ordinary}
If $X$ is not supersingular, there exists a quartic~$X_0$ defined over a
field $\Bbbk_0$, $\fchar\Bbbk_0=0$, with the property that
$\NS(X_0)\cong\NS(X)$, see \autoref{th.char=0}. Then, according to~\cite{DIS},
either $\ls|\Fn X_0|\le60$ or $\NS(X_0)=\bX_{64}$, and the same dichotomy
applies to the original surface~$X$.
If $\kchar=2$ or $3$, the last possibility is ruled out by
\autoref{cor.Schur}.
\qed

\remark\label{rem.ordinary}
According to~\cite{DIS} (and \autoref{th.char=0}), the number
$\ls|\Fn X|$ of lines in a
quartic~$X$ that is not supersingular takes values $\le52$, $54$, $56$, $60$,
or~$64$. Found in~\cite{DIS} is also a complete list of all configurations
$\Fn X$ of size at least~$54$.

Assume that $\kchar=2$ or~$3$.
Then, the maximal value $\ls|\Fn X|=64$ is ruled out by \autoref{th.ordinary}.
The next value $\ls|\Fn X|=60$ is realized by two configurations, $\bX_{60}'$
and $\bX_{60}''$ in the notation of~\cite{DIS}.
An explicit defining equation of the latter quartic is obtained in
\cite{rams.schuett:char2}, and it has a nonsingular reduction
(still with $60$ lines) over~$\F_4$.
In characteristic~$3$, the \emph{known} quartics become singular.
Conjecturally (M.~Sch\"{u}tt, private communication; the conjecture is based on
the arithmetical properties of the discriminant $\det\NS(X)=-60$ or $-55$),
even if a quartic~$X$
with $\Fn X\cong\bX_{60}'$ or $\bX_{60}''$
admits a nonsingular reduction modulo~$3$, the latter must be
supersingular;
by \autoref{th.char=3}, it would be isomorphic to the Fermat quartic.

One of the three configurations with $56$ lines is $\bX_{56}$; as was
recently
observed by T.~Shioda, over~$\C$ this surface is a non-standard projective
model of the Fermat quartic.
The defining equation
found by I.~Shimada~\cite{Shimada:X56}
has a nonsingular reduction
in characteristic~$3$
and the quartic obtained has $56$ lines.
Thus, the only case still open is that of
$\kchar=3$ and $\ls|\Fn X|=60$.
\endremark

\subsection{Quartics with a pencil of type $(10,0)$}\label{s.eq.(10,0)}
We conclude this section with the defining equations of a few supersingular
quartics in characteristic~$3$.

Any quartic $X\subset\Cp3$ containing a pencil $\pencil(b_2)$ 
of type $(10,0)$ is
supersingular and one has $\kchar=3$, see \autoref{prop.Euler};
arithmetically, such quartics are described in \autoref{prop.(10,0)}.
If $\pencil(b_2)$ has a section $b_1\in\Fn X$, this section intersects ten
lines $a_1,\ldots,a_{10}\in\pencil(b_2)$. Hence, $X$ satisfies the hypotheses of
\autoref{prop.10-2} and its equation is given by~\eqref{eq.10-2} as
\[
z_1^3z_3f\Bigl(\frac{z_0}{z_1},\frac{z_2}{z_3}\Bigr)
=z_0z_3^3-z_1z_2^3,
\label{eq.(10.0)}
\]
where $f(x,y)$ is an irreducible polynomial of bidegree~$(3,1)$ such that all
ten
roots of $f(y^3,y)$ are simple.
(Equations of this form have been studied in
\cite{Shimada:families}.)
As an immediate consequence,
we have the following
statement. (Note that, according to D.~Veniani, private communication, the conclusion of
this statement holds without the assumption that the pencil has a section.)

\corollary\label{cor.4-point}
If a pencil $\pencil(b_2)$ of type $(10,0)$ has a section, then, for any
$3$-fiber $s_1,s_2,s_3$ of the pencil, the four lines $b_1,s_1,s_2,s_3$
intersect at a single point.
\endcorollary

The projections of the ten common lines of the pencils
$\pencil(b_1)$, $\pencil(b_2)$ are the points
$A_i:=\pi(a_i)$ with coordinates $(y_i^3,y_i)$, $i=1,\ldots,10$,
where $y_i$ are the ten roots of the polynomial
$f(y^3,y)$. If $\kchar=3$,
one has
$(y_1^3,y_2^3;y_3^3,y_4^3)=(y_1,y_2;y_3,y_4)^3$;
hence,
\autoref{prop.10-2.lines} establishes a two-to-one
correspondence between the lines $l\in\Fn X$ disjoint from $b_1,b_2$ and the
quadruples $y_1,\ldots,y_4$ of roots of $f(y^3,y)$ satisfying the
equation $(y_1,y_2;y_3,y_4)=-1$.
Observe that (still assuming $\kchar=3$)
\roster*
\item
one has $(y_1,y_2;y_3,\infty)=-1$ if and only if $y_1+y_2+y_3=0$, and
\item
if $y_i',y_i''$ are the roots of $y^2+p_iy+q_i$, $i=1,2$, then
 $(y_1',y_1'';y_2',y_2'')=-1$ if and only if $p_1p_2+q_1+q_2=0$.
\endroster
Below, we consider a few special cases.

\proposition\label{prop.(10,0)+(10,0)}
Any quartic $X$ containing a pair of obverse pencils of type $(10,0)$ is
projectively equivalent to the Fermat quartic.
\endproposition

\proof
Since $\pencil(b_1)$ is also of type~$(10,0)$, one has $f(x,y)=x^3-y$
in~\eqref{eq.(10.0)}.
\endproof

\proposition\label{prop.(10,0)+(4,6)}
A quartic~$X$ is as in \autoref{th.char=3}\iref{58.(10,0)} if and only if
$X$ contains a pair of obverse pencils $\pencil(b_1)$, $\pencil(b_2)$ of
types $(4,6)$ and $(10,0)$ respectively.
Up to projective transformation, such quartics constitute the
$1$-parameter family given by \eqref{eq.(10.0)} with
\[*
f(x,y)=
x^3-(\pv+1)(\pv y+y+\pv)x^2+\pv^2(\pv+1)(y+\pv+1)x-\pv^4y,
\]
where $\pv\in\Bbbk\sminus\F_3$.
\endproposition

\proof
The ``only if'' part is given by the explicit description of the
configuration of lines in~$X$. For the converse, we use the computation in
the proof of \autoref{lem.Schur.uniqueness}. In the notation introduced
there, let $B_2$ be given by $y^3=x$; then, in appropriate coordinates,
$B_1$ is in the $1$-parameter family $B_1(u)$, $u\ne1$, considered in the
proof.
To simplify the notation, we let $u:=\pv-1$; then $\pv\ne-1$ (as otherwise
$B_1$ is purely inseparable and $X$ is the Fermat quartic, see
\autoref{prop.(10,0)+(10,0)}) and $\pv\ne0,1$ (as otherwise $B_1$ is
reducible). By the construction, the intersection $B_1\cap B_2$ contains the
points $(r,r)$, $r=0,1,\infty$, and the fourth point $(\Gm,\Gl)=(\pv^3,\pv)$
found in the proof also lies in~$B_2$. (In other words, assuming that
$\pencil(b_1)$ has three $3$-fibers, we obtain a fourth one, which agrees
with \autoref{cor.exotic.type}.)

In addition to $0$, $1$, $\pv$, $\infty$, the polynomial $f(y^3,y)$ has three
pairs of roots $(y_i',y_i'')$ satisfying the quadratic equations
$q_i(y_i')=g_i(y_i'')=0$, where $i=0,1,\pv$ and
\[*
g_0(y):=y^2-\pv,\qquad
g_1(y):=y^2+y+\pv,\qquad
g_\pv(y):=y^2+\pv y+\pv.
\]
Using the two observations prior to \autoref{prop.(10,0)+(10,0)}, one can see that
\[*
(i,\infty;y_i',y_i'')=(i,j;y_k',y_k'')=(y_i',y_i'';y_k',y_k'')=-1
\]
for any permutation $\{i,j,k\}$ of $\{0,1,\pv\}$.
It follows that $X$ contains at least $58$ lines and, thus, is as in
\autoref{th.char=3}\iref{58.(10,0)}.
\endproof

\proposition\label{prop.(10,0)++}
Up to projective equivalence, the quartics as in
\autoref{th.char=3}\iref{58.(10,0)x2}
constitute the $1$-parameter family given by~\eqref{eq.(10.0)} with
\[*
f(x,y)=\pv y - (\pv+1)x + x^3,
\]
where $\pv\in\Bbbk^\times$.
\endproposition

\proof
Let $X$ be a quartic as in the statement; then $\Fn X$ is the union of the
two type $(10,0)$ pencils contained in~$X$. Choose for $\pencil(b_2)$ one of
these pencils; then, the other pencil is $\pencil(a_\infty)$ for a certain line
$a_\infty\in\pencil(b_2)$, and each section of $\pencil(b_2)$
intersects~$a_\infty$. Pick a section~$b_1$ and consider the corresponding
projection $\pi:=\pi[b_1,b_2]\:X\to\Cp1\times\Cp1$. The coordinates $(x,y)$
in $\Cp1\times\Cp1$ can be chosen so that
$B_2:=\pi(b_2)$ is given by $x=y^3$,
the line~$a_\infty$ projects to
$A_\infty(\infty,\infty)$, and some other line $a_0\in\pencil(b_1)\cap\pencil(b_2)$
projects to $(0,0)$.

According to \autoref{cor.4-point}, the curve $B_1:=\pi(b_1)$ is inflection
tangent at $A_\infty$ to the fiber $y=\infty$; hence, its
defining polynomial is of the form
\[*
f(x,y)=\pv y+ux+vx^2+x^3,\quad \pv\ne0.
\]
A simple count shows that there is a line $a\in\pencil(b_1)\cap\pencil(b_2)$
other than~$a_\infty$ that intersects at least
eight lines disjoint from
$b_1,b_2$.
Assuming that $a=a_0$,
from the observation prior to
\autoref{prop.(10,0)+(10,0)} one concludes that the polynomial $f(y^3,y)$
must be odd; hence, $v=0$. Up to projective transformation, we can also
assume that $\pm1$ are among the roots; then $u=-(\pv+1)$ and $f(x,y)$ is as
in the statement.

In addition to $0$, $\infty$, and $\pm1$, the roots of $f(y^3,y)$ are those
of $y^6+y^4+y^2-\pv$;
denoting by~$t$ one of the extra roots, one can easily see that all six roots are
$\pm t$ and $\pm(t\pm1)$; they are all pairwise distinct and
different from $0$ and $\pm1$. (Since $\pv\ne0$, one has $t\notin\F_3$).
The set $\{\pm1,\pm t,\pm(t\pm1)\}$ contains eight triples summing up to~$0$.
Hence, $X$ contains at least $58$ lines and, thus, is as in the
statement.
\endproof

{
\let\.\DOTaccent
\def\cprime{$'$}
\bibliographystyle{amsplain}
\bibliography{degt}

\providecommand{\bysame}{\leavevmode\hbox to3em{\hrulefill}\thinspace}
\providecommand{\MR}{\relax\ifhmode\unskip\space\fi MR }
\providecommand{\MRhref}[2]{%
  \href{http://www.ams.org/mathscinet-getitem?mr=#1}{#2}
}
\providecommand{\href}[2]{#2}
\begin{thebibliography}{10}

\bibitem{Artin:supersingular}
M.~Artin, \emph{Supersingular {$K3$} surfaces}, Ann. Sci. \'Ecole Norm. Sup.
  (4) \textbf{7} (1974), 543--567 (1975). \MR{0371899 (51 \#8116)}

\bibitem{Barth:K3}
W.~Barth, \emph{Lectures on {$K3$}- and {E}nriques surfaces}, Algebraic
  geometry, {S}itges ({B}arcelona), 1983, Lecture Notes in Math., vol. 1124,
  Springer, Berlin, 1985, pp.~21--57. \MR{805328 (86m:14027)}

\bibitem{Bourbaki:Lie}
Nicolas Bourbaki, \emph{Lie groups and {L}ie algebras. {C}hapters 4--6},
  Elements of Mathematics (Berlin), Springer-Verlag, Berlin, 2002, Translated
  from the 1968 French original by Andrew Pressley. \MR{1890629 (2003a:17001)}

\bibitem{Conway.Sloane}
J.~H. Conway and N.~J.~A. Sloane, \emph{Sphere packings, lattices and groups},
  Grundlehren der Mathematischen Wissenschaften [Fundamental Principles of
  Mathematical Sciences], vol. 290, Springer-Verlag, New York, 1988, With
  contributions by E. Bannai, J. Leech, S. P. Norton, A. M. Odlyzko, R. A.
  Parker, L. Queen and B. B. Venkov. \MR{920369 (89a:11067)}

\bibitem{DIK}
Alex Degtyarev, Ilia Itenberg, and Viatcheslav Kharlamov, \emph{Real {E}nriques
  surfaces}, Lecture Notes in Mathematics, vol. 1746, Springer-Verlag, Berlin,
  2000. \MR{1795406 (2001k:14100)}

\bibitem{DIS}
Alex Degtyarev, Ilia Itenberg, and Ali~Sinan Sert\"oz, \emph{Lines on
  quartics}, to appear, \verb+arXiv:1601.04238+, 2016.

\bibitem{GAP4}
\emph{{GAP} {\textendash} {G}roups, {A}lgorithms, and {P}rogramming, {V}ersion
  4.7.7}, \href {http://www.gap-system.org}
  {\texttt{http://www.gap-system.org}}, Feb 2015.

\bibitem{Huybrechts}
Daniel Huybrechts, \emph{Lectures on {$K3$} surfaces}, Cambridge University
  Press, to appear.

\bibitem{Kulikov:periods}
Vik.~S. Kulikov, \emph{Surjectivity of the period mapping for {$K3$} surfaces},
  Uspehi Mat. Nauk \textbf{32} (1977), no.~4(196), 257--258. \MR{0480528 (58
  \#688)}

\bibitem{Lieblich.Maulik}
Max Lieblich and Davesh Maulik, \emph{A note on the cone conjecture for ${K}3$
  surfaces in positive characteristic}, \verb+arXiv:1102.3377+.

\bibitem{Miranda.Morrison:book}
Rick Miranda and David~R. Morrison, \emph{Embeddings of integral quadratic
  forms}, electronic, \url{http://www.math.ucsb.edu/~drm/manuscripts/eiqf.pdf},
  2009.

\bibitem{Nikulin:forms}
V.~V. Nikulin, \emph{Integer symmetric bilinear forms and some of their
  geometric applications}, Izv. Akad. Nauk SSSR Ser. Mat. \textbf{43} (1979),
  no.~1, 111--177, 238, English translation: Math USSR-Izv. 14 (1979), no. 1,
  103--167 (1980). \MR{525944 (80j:10031)}

\bibitem{Payne.Thas}
Stanley~E. Payne and Joseph~A. Thas, \emph{Finite generalized quadrangles},
  second ed., EMS Series of Lectures in Mathematics, European Mathematical
  Society (EMS), Z\"urich, 2009. \MR{2508121 (2010k:51013)}

\bibitem{rams.schuett:char2}
S{\l}awomir Rams and Matthias Sch{\"u}tt, \emph{At most 64 lines on smooth
  quartic surfaces (characteristic 2)}, to appear, \verb+arXiv:1512.01358+,
  2012.

\bibitem{rams.schuett:char3}
\bysame, \emph{112 lines on smooth quartic surfaces (characteristic 3)}, Q. J.
  Math. \textbf{66} (2015), no.~3, 941--951. \MR{3396099}

\bibitem{rams.schuett}
\bysame, \emph{64 lines on smooth quartic surfaces}, Math. Ann. \textbf{362}
  (2015), no.~1-2, 679--698. \MR{3343894}

\bibitem{Rudakov.Shafarevich}
A.~N. Rudakov and I.~R. Shafarevich, \emph{Surfaces of type {$K3$} over fields
  of finite characteristic}, Current problems in mathematics, {V}ol. 18, Akad.
  Nauk SSSR, Vsesoyuz. Inst. Nauchn. i Tekhn. Informatsii, Moscow, 1981,
  pp.~115--207. \MR{633161 (83c:14027)}

\bibitem{Saint-Donat}
B.~Saint-Donat, \emph{Projective models of {$K$-$3$} surfaces}, Amer. J. Math.
  \textbf{96} (1974), 602--639. \MR{0364263 (51 \#518)}

\bibitem{Schur:quartics}
Friedrich Schur, \emph{Ueber eine besondre {C}lasse von {F}l\"achen vierter
  {O}rdnung}, Math. Ann. \textbf{20} (1882), no.~2, 254--296. \MR{1510168}

\bibitem{Segre}
B.~Segre, \emph{The maximum number of lines lying on a quartic surface}, Quart.
  J. Math., Oxford Ser. \textbf{14} (1943), 86--96. \MR{0010431 (6,16g)}

\bibitem{Shimada:families}
Ichiro Shimada, \emph{On supercuspidal families of curves on a surface in
  positive characteristic}, Math. Ann. \textbf{292} (1992), no.~4, 645--669.
  \MR{1157319}

\bibitem{Shimada:X56}
Ichiro Shimada and Tetsuji Shioda, \emph{On a smooth quartic surface containing
  56 lines which is isomorphic as a {$K3$} surface to the {F}ermat quartic.},
  to appear, 2016.

\bibitem{vanderBlij}
F.~van~der Blij, \emph{An invariant of quadratic forms mod {$8$}}, Nederl.
  Akad. Wetensch. Proc. Ser. A 62 = Indag. Math. \textbf{21} (1959), 291--293.
  \MR{0108467 (21 \#7183)}

\end{thebibliography}
}

\end{document}